\documentclass[11pt,a4paper,reqno]{amsart}

\pdfoutput=1

\usepackage{enumerate}
\usepackage[nocompress]{cite}
\usepackage{xcolor}
\usepackage{hyperref}
\hypersetup{
    colorlinks,
    linkcolor={red!50!black},
    citecolor={blue!50!black},
    urlcolor={blue!80!black},
    linktoc=page}
\usepackage{tikz-cd}
\usepackage{amssymb}
\usepackage{graphicx}
\usepackage{amscd}
\usepackage[margin=3.14cm]{geometry}
\usepackage{mathtools}
\usepackage{bm, bbm}
\usepackage{aliascnt}
\usepackage[symbol,perpage]{footmisc}
\usepackage[capitalise]{cleveref}

\usepackage{booktabs}
\usepackage{array}

\usepackage[
    final,
    kerning=true,
    spacing=true,
    factor=1100,
    stretch=10,
    shrink=10]{microtype}

\theoremstyle{plain}
\newtheorem{theorem}{Theorem}[section]
\numberwithin{equation}{section}
\numberwithin{figure}{section}
\numberwithin{table}{section}
\newtheorem*{theorem*}{Theorem}
\newtheorem{maintheorem}{Theorem}

\crefname{maintheorem}{Theorem}{Theorems}
\Crefname{maintheorem}{Theorem}{Theorems}
\newaliascnt{proposition}{theorem}
\newtheorem{proposition}[proposition]{Proposition}
\aliascntresetthe{proposition}
\newaliascnt{lemma}{theorem}
\newtheorem{lemma}[lemma]{Lemma}
\aliascntresetthe{lemma}
\newaliascnt{corollary}{theorem}
\newtheorem{corollary}[corollary]{Corollary}
\aliascntresetthe{corollary}
\newaliascnt{conjecture}{theorem}

\aliascntresetthe{conjecture}
\newaliascnt{problem}{theorem}

\aliascntresetthe{problem}
\theoremstyle{definition}
\newaliascnt{definition}{theorem}
\newtheorem{definition}[definition]{Definition}
\aliascntresetthe{definition}
\newtheorem*{definition*}{Definition}
\newaliascnt{remark}{theorem}
\newtheorem{remark}[remark]{Remark}
\aliascntresetthe{remark}
\newaliascnt{example}{theorem}
\newtheorem{example}[example]{Example}
\aliascntresetthe{example}
\newaliascnt{assumption}{theorem}

\aliascntresetthe{assumption}

\newcommand		{\p}[1]	{\left(#1\right)}
\newcommand		{\abs}[1]{\left|#1\right|}
\newcommand 	{\norm}[1]{\left\lVert#1\right\rVert}
\newcommand    	{\R} {\mathbb R}
\newcommand    	{\N} {\mathbb N}
\newcommand    	{\Q} {\mathbb Q}
\newcommand    	{\Z} {\mathbb Z}
\newcommand    	{\T} {\mathbb T}
\DeclareMathOperator* {\E} {\mathbb E}
\DeclareMathOperator {\var} {\mathrm {Var}}
\newcommand    	{\ve} {\varepsilon}
\renewcommand    	{\d} {\, \mathrm d}
\newcommand    	{\1} {\mathbbm 1}
\renewcommand{\emph}[1]{\textit{#1}}
\newcommand    	{\id} {\mathrm{id}}
\newcommand    	{\lip} [1] {\mathrm{Lip}\p{#1}}

\newcommand 	{\CM} {\mathcal M}
\renewcommand 	{\CD} {\mathcal D}
\newcommand 	{\CO} {\mathcal O}
\newcommand 	{\CL} {\mathcal L}
\newcommand 	{\CT} {\mathcal T}
\newcommand 	{\CW} {\mathcal W}

\newcommand 	{\clip} {\mathrm{clip}}
\DeclareMathOperator {\dist} {dist}

\newcommand{\Linf}[1]{L^\infty_{[-1,1]}(#1)}

\title[Nonlinear action convergence]
{Action convergence of nonlinear operators and coupled dynamics
on graph limits}

\author{Davide Sclosa, Matteo Tanzi}

\begin{document}

\begin{abstract}
We extend to nonlinear maps the notion of action convergence introduced by
Backhausz and Szegedy for linear maps.
Action convergence provides a common framework for limits
of dense, sparse and intermediate graph sequences. We prove that an action convergent
sequence of nonlinear transformations of the~$L^\infty$ unit ball has a limit if and
only if it is $L^1$-equicontinuous. For coupled maps on graphs, we show that action
convergence of the linear graph operators implies action convergence of the
nonlinear dynamics, with no assumption on the initial conditions; the same holds
for coupled ODEs. We then derive consequences for large finite networks:
the self-consistent transfer
operator is asymptotically exact on quasirandom graphs, and macroscopic Lyapunov
contraction in the limit yields approximate synchronization of the finite systems, with
an error uniform in time. Moreover, for any
nondegenerate choice of the coupling functions, the limit of the nonlinear coupled maps
determines the limit of the graphs.
We also describe explicit limits for cycles, recursive scale-free trees and quasirandom graphs.
\end{abstract}

\vspace*{-0.1cm}
\maketitle


\vspace{-0.5cm}
\setcounter{tocdepth}{1}
\tableofcontents

\vspace{-0.6cm}
\section{Introduction}

Graph limit theory compares graphs with different numbers of vertices by representing
their limits as analytic objects. It has led to advances in extremal graph theory, large
deviations for random graphs and statistical inference for networks, and has more recently
been used to compare graph neural networks with their limiting models.
Separate theories were initially developed for dense graph
sequences~\cite{lovasz2012large} and for sparse ones~\cite{hatami2014limits}. In a recent
seminal paper, Backhausz and Szegedy~\cite{backhausz2022action} introduced a unified theory
that treats dense, sparse and intermediate densities simultaneously.
Their viewpoint is both analytic and dynamical: each graph on~$n$ vertices is associated
with a linear operator~$W_n$, such as its adjacency matrix or the transition matrix of its
random walk, and the graph is studied through the \emph{action} of~$W_n$ on vertex states,
that is, through the linear dynamics
\begin{equation}\label{eq:intro_lin}
	u(i,t+1)=\sum_{j=1}^n W_n(i,j)\,u(j,t),
\end{equation}
where~$u(i,t)\in\R$ is the state of vertex~$i\in\{1,\ldots,n\}$ at time~$t\geq0$.
A graph sequence converges if its actions converge as~$n\to\infty$, in a
suitable sense called \emph{action convergence}.

As Backhausz and Szegedy already observed, action convergence does not require linearity.
It has since been extended to multilinear maps~\cite{zucal2023action}.
In this paper, we extend the theory to nonlinear
maps. We prove that an action convergent sequence of nonlinear maps has a limit if and
only if it is equicontinuous~(\cref{thm:existence}), that action convergence is preserved
under continuous pairwise coupling~(\cref{thm:pairwise}), and that convergence of the
one-step maps implies joint convergence over every finite time
horizon~(\cref{thm:trajectories}).
Nonlinear maps were independently considered in the context of graph neural networks, but the
existence of action limits was left open~\cite[Section~4.4]{le2023limits}.

Our results extend action convergence from graph operators to nonlinear maps built
from them, including graph neural network layers and nonlinear coupled systems.
Our main motivation is the latter, in discrete time, although our results also apply to
coupled ODEs (\cref{sec:continuous_time}).
For continuous maps~$F$ and~$G$, we consider
\begin{equation}\label{eq:intro_nonlin}
	u(i,t+1)
	=F\Bigl(u(i,t),\sum_{j=1}^n W_n(i,j)\,
	G\bigl(u(i,t),u(j,t)\bigr)\Bigr).
\end{equation}
Under the equicontinuity assumption of~\cref{thm:pairwise}, action convergence of the
graphs, that is of~\eqref{eq:intro_lin}, implies action convergence
of~\eqref{eq:intro_nonlin}.
Conversely, under mild conditions on the graph operators and on~$F,G$, the limit
of~\eqref{eq:intro_nonlin} determines the limit of the graphs (\cref{prop:converse}):
in other words, one can \textit{``hear" the shape of a network} through any nondegenerate
dynamics on it.
We illustrate the usefulness of action convergence by
deriving several consequences for~\eqref{eq:intro_nonlin} on large graphs,
with no assumptions on the initial condition:
the self-consistent transfer operator, which describes the
evolution of the distribution of states and is exact on complete graphs, is asymptotically
exact on quasirandom graphs~(\cref{cor:STO}). Laws of approximate fixed points pass to and
from the limit~(\cref{cor:fixed}). Macroscopic Lyapunov estimates in the limit give a
criterion for approximate synchronization on large finite graphs, uniformly in
time~(\cref{cor:synch}).


\subsection{Main idea}\label{subsec:intro_idea}
We now introduce action convergence informally, through the lens
of~\eqref{eq:intro_nonlin}, and explain its advantages.
We consider vertex states in~$[-1,1]$ and continuous
functions~$F,G\colon[-1,1]^2\to[-1,1]$, although this can be generalized
(\cref{rem:compact_states}). The dynamics~\eqref{eq:intro_nonlin} is then given by
iterating the map~$T_n\colon[-1,1]^n\to[-1,1]^n$ defined by
\begin{equation}\label{eq:network}
	(T_nu)(i):=F\Bigl(u(i),\ \sum_{j=1}^nW_n(i,j)\,G\bigl(u(i),u(j)\bigr)\Bigr),
	\qquad 1\leq i\leq n.
\end{equation}
Since the maps~$(T_n)_{n\in \N}$ act on different spaces, there is no immediate way to
define convergence.
Fix a state~$u=(u(i))_{1\leq i\leq n} \in [-1,1]^n$. The empirical
measure~$\CD(u,T_nu):= \frac{1}{n} \sum_{i=1}^n\delta_{(u(i),(T_nu)(i))}$
is the law of the input-output~$(u(i),(T_nu)(i))$ with
vertex~$i$ sampled uniformly from~$\{1,\ldots,n\}$.
It is a probability measure on~$[-1,1]^2$ regardless
of~$n$, and this is what generalizes.
For a general probability space~$(\Omega,\mu)$, playing the role of a set of vertices,
states are measurable functions~$w\colon \Omega\to [-1,1]$,
and dynamics is a transformation~$T$ mapping a state~$w$ to a state~$Tw$.
Let~$\CD(w,Tw)$ be the law of~$(w(x),Tw(x))$ with~$x\sim \mu$,
i.e., the push-forward of~$\mu$ through~$x\mapsto (w(x),Tw(x))$.
It is again a probability measure on~$[-1,1]^2$.
Roughly speaking, we call~$T$ a limit of~$(T_n)$ if for every~$\ve>0$, the following is true for all~$n$ sufficiently
large:
for \textit{every} state~$u\in[-1,1]^n$ of the finite system,
there is a state~$w\in\Linf{\mu}$ of the limit system such
that~$\CD_{\mu_n}(u,T_nu)$ and~$\CD_{\mu}(w,Tw)$ are at most $\ve$-far apart
in the Wasserstein metric, and conversely for \textit{every}
state~$w\in\Linf{\mu}$ of the limit system there is a state~$u\in[-1,1]^n$ of the finite system
such that~$\CD_{\mu_n}(u,T_nu)$ and~$\CD_{\mu}(w,Tw)$ are at most $\ve$-far apart.
The first approximation says that every behaviour of the large finite systems is
reproduced in the limit; the second, that the limit adds nothing spurious, so that
whatever we find in the limit is meaningful for the finite system.
In the actual definition,
the law of $k$ states are considered at once,
to record their correlations (\cref{def:profile}),
and vertices~$i\in\{1,\ldots,n\}$ can be sampled non-uniformly, e.g.
giving a hub weight that does not vanish
as~$n\to\infty$ (\cref{subsec:scale_free}).

Large-$n$ descriptions of coupled systems such as~\eqref{eq:intro_nonlin} have been the object of intense recent investigation. Compared with previous approaches,
ours has three advantages.

\emph{(1) No assumptions on the initial condition.}
	Previous approaches follow, for each~$n$, a single trajectory,
	and compare it with a limit
equation (Vlasov, McKean--Vlasov, or a self-consistent transfer operator) under strong
assumptions on the initial condition.
The initial condition is assumed random and independent across
vertices~\cite{delattre2016note,jabin2025mean,maquera2025selfconsistent},
which can make the structure of the graph invisible~\cite[Section~1.2]{coppini2020law},
or it is the
discretization of a continuum
initial condition~\cite{medvedev2014nonlinearS,medvedev2014nonlinear,
kuehn2020network,gkogkas2022graphop,bick2025dynamical},
which requires an ordering of the vertices under which both the initial
condition and the graph operator converge in norm simultaneously.
In contrast, action convergence considers all initial conditions for every~$n$.

\emph{(2) Robustness to small changes of the graph.}
	By design,
	action convergence only retains the robust,
	macroscopic features of a system and discards
	the fragile, microscopic ones.
	For example, adding or removing a single edge can destroy a fixed point
	of~$T_n$, no matter how large the graph is. In contrast, approximate
	fixed points are robust under
	such perturbations (\cref{subsec:fixed}).

\emph{(3) All densities at once.}
Action convergence covers dense, sparse and intermediate densities alike.
For dense graphs, it is equivalent to graphon convergence in cut
norm; for graphs of bounded degree, it coincides with local--global
convergence~\cite{hatami2014limits,backhausz2022action}, a strengthening of the local
weak convergence
used for interacting dynamics in~\cite{oliveira2020interacting,lacker2023local}.


\subsection{Main results}\label{subsec:definitions}

Our main theorems are stated and proved for general maps~$T$ acting on measurable
functions~$u$, which need not come from graphs. This analytic point of view makes the
key ideas transparent. Nonlinear maps arising from graphs remain, however, the
motivation and main application of our theory, in the same way as linear maps arising
from graphs are the main application in~\cite{backhausz2022action}.

Let $(\Omega,\mu)$ denote a standard Borel probability space,
such as the Lebesgue unit
interval~$([0,1],\lambda)$, or the finite
spaces~$([n],\mu_n)$, where~$[n]:=\{1,\ldots,n\}$ and~$\mu_n(i):=\mu_n(\{i\})>0$.
We write~$\Linf{\mu}$ for the set
of measurable maps~$u\colon\Omega\to[-1,1]$, identified up to $\mu$-null sets.
On~$([n],\mu_n)$ this is the cube~$[-1,1]^n$,
and we write~$u=(u(1),\ldots,u(n)) \in [-1,1]^n$.
We denote by~$\CT$ the class of triples
\[
\CT :=
\left\{
(\Omega,\mu,T)\ \middle|\
\begin{array}{l}
(\Omega,\mu)\text{ is a standard Borel probability space},\\
T\colon\Linf{\mu}\to\Linf{\mu}\text{ is }
\text{$L^1(\mu)$-uniformly continuous}
\end{array}
\right\}.
\]
We often write~$T\in\CT$, leaving the underlying space implicit.
If~$T\in\CT$,
then the iterates~$T^2,T^3,\ldots$ belong to~$\CT$, so
every~$T\in\CT$ generates dynamics on~$\Linf{\mu}$.

We denote by~$\CL\subseteq\CT$ the subclass of \emph{linear} maps, that is, of
the~$(\Omega,\mu,W)\in\CT$ such that~$W$ is the restriction to~$\Linf{\mu}$ of a linear
operator on~$L^\infty(\mu)$. Equivalently, $W$ is a linear operator on~$L^\infty(\mu)$
with~$\norm{Wu}_\infty\leq\norm{u}_\infty$, such that~$\norm{Wu_k}_{L^1(\mu)}\to0$
whenever~$u_k\in\Linf{\mu}$ and~$\norm{u_k}_{L^1(\mu)}\to0$.
Let~$W_n$ be a real matrix with~$\sum_j\abs{W_n(i,j)}\leq1$ for every~$i$,
so that~$[-1,1]^n$ is invariant. Then~$([n],\mu_n,W_n) \in \CL$; the matrix
need not be nonnegative, nor symmetric. The map~$T_n$
in~\eqref{eq:network} satisfies~$([n],\mu_n,T_n) \in \CT$.
Note that~$W_n$ and~$T_n$
generate the dynamics~\eqref{eq:intro_lin} and~\eqref{eq:intro_nonlin}, respectively.

Here and below, $Tu(x)$ stands for~$(Tu)(x)$. For~$k\geq1$
and~$u=(u_1,\ldots,u_k)\in\Linf{\mu}^k$, we denote by~$\CD_\mu(u)$ the law
of~$u$, namely the push-forward of~$\mu$ under the map~$x\mapsto(u_1(x),\ldots,u_k(x))$.
We let~$T$ act componentwise, $Tu:=(Tu_1,\ldots,Tu_k)$, and write~$\CD_\mu(u,Tu)$
for the law of~$(u,Tu)$. It is an element of~$\CM([-1,1]^{2k})$, the
space of Borel probability measures on~$[-1,1]^{2k}$. For~$([n],\mu_n)$,
it is~$\sum_{i=1}^n\mu_n(i)\,\delta_{(u(i),Tu(i))}$.

\begin{definition}\label{def:profile}
For~$k\geq1$, the \emph{$k$-profile} of a map~$T\colon\Linf{\mu}\to\Linf{\mu}$ is
\[
	P_k(T):=\{\CD_\mu(u,Tu)\colon u\in\Linf{\mu}^k\}.
\]
\end{definition}

In words, profiles record the input--output statistics of~$T$ for all inputs.
The $1$-profile does so for single inputs, the $k$-profile for $k$ inputs
simultaneously.

Fix~$m\geq 1$. The
Wasserstein metric~$d_W(\rho,\nu):=\sup_h \int h \d(\rho-\nu)$,
where~$h$ varies among functions
satisfying~$\abs{h(x)-h(y)}\leq \sum_{i=1}^{m} \abs{x_i-y_i}$,
turns~$\CM([-1,1]^{m})$ into a compact metric space.
The corresponding Hausdorff pseudometric
on subsets~$P,Q\subseteq \CM([-1,1]^{m})$
is given by~$d_H(P,Q):=\max\{\sup_{\rho\in P}d_W(\rho,Q),\sup_{\rho\in Q}d_W(\rho,P)\}$.
Restricted to compact subsets, the pseudometric~$d_H$ is a genuine metric,
and the resulting topology is compact.
Let~$\overline P$ denote the closure in the Wasserstein metric.
We have~$d_H(P,Q)=d_H(\overline{P},\overline{Q})$.

\begin{definition}\label{def:action_convergence}
A sequence~$(T_n)_{n\in\N}$ in~$\CT$ is \emph{action convergent} if, for
every~$k\geq1$, the sets~$\overline{P_k(T_n)}$ converge in the Hausdorff metric.
A map~$T\in\CT$ is a \emph{limit} of the sequence, and we write~$T_n\to T$,
if~$\overline{P_k(T_n)}\to\overline{P_k(T)}$ for every~$k\geq1$.
\end{definition}

We call~$(\Omega,\mu,T)$ and~$(\Omega',\mu',T')$ \emph{weakly equivalent}
if~$\overline{P_k(T)}=\overline{P_k(T')}$ for all~$k\geq1$. Weakly equivalent maps
are indistinguishable in this theory, and limits are unique exactly up to weak
equivalence. If~$T_n\to T$, then~$T_n\to T'$ if and
only if~$T'$ is weakly equivalent to~$T$. ``The limit'' of a sequence
always means a representative of this class.

If~$(T_n)$ is action convergent, the sets~$\overline{P_k(T_n)}$ have limits, but these
need not be the profiles of a single map: an action convergent sequence need not have
a limit in~$\CT$. For linear operators, Backhausz and
Szegedy~\cite{backhausz2022action} proved that limits exist under a uniform bound
on~$\norm{\cdot}_{p\to q}$ for some~$p<\infty$.
Our first result is a nonlinear analogue with~$p=1$.
A sequence~$(T_n)$ in~$\CT$ is \emph{$L^1$-equicontinuous}
if there is a modulus~$\omega$ (\cref{def:modulus}) such
that~$\norm{T_nu-T_nw}_{L^1(\mu_n)}\leq\omega(\norm{u-w}_{L^1(\mu_n)})$ for all~$u,w\in\Linf{\mu_n}$
and all~$n$.

\begin{maintheorem} \label{thm:existence}
An action convergent sequence in~$\CT$ has a limit in~$\CT$ if and only if it is
$L^1$-equicontinuous. In that case, every limit has the same modulus of continuity.
\end{maintheorem}

A detailed comparison between \cref{thm:existence} and the theorem of Backhausz and Szegedy
is discussed in \cref{rem:linear_relation}, where we also sketch the proof of a
common generalization.

Our second result shows that action convergence of linear maps is preserved under
pairwise coupling with continuous functions.
Every~$W\in\CL$ can be represented as
\[
	Wu(x)=\int_\Omega u\d\nu_x
	\qquad\text{for $\mu$-almost every~$x$},
\]
for a measurable family~$(\nu_x)_{x\in\Omega}$ of finite signed measures
with~$\abs{\nu_x}(\Omega)\leq 1$ (\cref{lem:signed_kernel}).

\begin{maintheorem}[Pairwise coupling] \label{thm:pairwise}
Let~$F,G\colon[-1,1]^2\to[-1,1]$ be continuous. For~$W\in\CL$,
define the map~$\Gamma(W)$ by
\[
	(\Gamma(W)u)(x):=F\Bigl(u(x),\,\int_\Omega G\bigl(u(x),u(y)\bigr)\d\nu_x(y)\Bigr),
\]
where~$Wu(x)=\int_\Omega u\d\nu_x$. Then~$\Gamma(W)\in\CT$.
Moreover, if~$(W_n)$ is an $L^1$-equicontinuous action convergent sequence in~$\CL$,
then it has a limit~$W\in\CL$, and~$\Gamma(W_n)\to\Gamma(W)$.
\end{maintheorem}

If~$W_n$ is a matrix, then~$\nu_i=\sum_jW_n(i,j)\delta_j$, and~$\Gamma(W_n)$ is the
map~$T_n$ in~\eqref{eq:network}.
Our third result concerns trajectories.
For~$t,k\geq1$, the \emph{trajectory profile} of~$T$ is
\[
	P^{(t)}_k(T):=\bigl\{\CD_\mu(u,Tu,\ldots,T^tu)\colon u\in\Linf{\mu}^k\bigr\}
	\subseteq\CM\bigl([-1,1]^{k(t+1)}\bigr).
\]
For~$t=1$, this is the $k$-profile~$P_k(T)$.

\begin{maintheorem}[Trajectories] \label{thm:trajectories}
Let~$(T_n)$ be an $L^1$-equicontinuous sequence in~$\CT$, action convergent to~$T$. Then,
for all~$t,k\geq1$,
\[
	d_H\bigl(P^{(t)}_k(T_n),P^{(t)}_k(T)\bigr)\to0.
\]
\end{maintheorem}

Projecting the trajectory profile
to first and last marginals gives $\CD_\mu(u,T^tu)$,
action convergence of the maps implies action convergence of every fixed
iterate.

By \cref{lem:moduli}, a sequence that converges to a limit in~$\CT$ is
automatically $L^1$-equi\-con\-ti\-nuous; hence, in \cref{thm:pairwise,thm:trajectories},
equicontinuity can be replaced by the assumption that the limit lies in~$\CT$.
The meaning of \cref{thm:pairwise,thm:trajectories} for coupled dynamics
on graphs is explained in \cref{sec:applications}.

\subsubsection*{One can hear the shape of a network}
Under mild assumptions, the map~$W\mapsto\Gamma(W)$ is injective up to weak equivalence
(\cref{prop:converse}).
Assume that~$W,W'\in\CL$ satisfy~$W\1=W'\1=\1$, and that~$F,G$ satisfy a mild
nondegeneracy condition. If~$\Gamma(W)$ and~$\Gamma(W')$ are weakly equivalent, then so
are~$W$ and~$W'$. More precisely, if~$(W_n)$ is
$L^1$-equicontinuous
with~$W_n\1=\1$, then~$\Gamma(W_n)\to\Gamma(W)$ implies~$W_n\to W$
(\cref{prop:converse}). Combined with \cref{thm:pairwise}, convergence of~$\Gamma(W_n)$
for a single nondegenerate pair~$F,G$ implies its convergence for all
continuous~$F$ and~$G$.

\subsubsection*{Overview of the proofs}
\Cref{thm:pairwise}, its converse, and \cref{thm:trajectories},
all follow from a single ``lift and repair'' lemma, which is the main technical contribution
of this paper. We illustrate it here on \cref{thm:trajectories} with~$t=2$ and~$k=1$.
Trajectory laws~$\CD_\mu(u,Tu,T^2u)$ are not part of the profiles of~$T$, but they can be
read off from the $2$-profile: they correspond to the laws~$\CD_\mu(u_1,u_2,Tu_1,Tu_2)$
with~$u_2=Tu_1$, that is, laws in~$P_2(T)$ supported on the
set~$K:=\{(x_1,x_2,y_1,y_2)\in [-1,1]^4 \colon x_2=y_1\}$. Thus, it suffices to show that
convergence of the $2$-profiles passes to their intersections with~$\CM(K)$.
This is not automatic: a law~$\CD_\mu(u_1,u_2,Tu_1,Tu_2)$ of~$T$ with~$u_2=Tu_1$ is
approximated by laws~$\CD_{\mu_n}(w_1,w_2,T_nw_1,T_nw_2)$ of~$T_n$
with~$w_2\approx T_nw_1$: they satisfy the constraint only approximately. We
\emph{repair} them: replacing~$w_2$ by~$T_nw_1$ moves~$w_2$
by~$\norm{w_2-T_nw_1}_{L^1(\mu_n)}$, which is small, and hence moves~$T_nw_2$ by a small
amount, by equicontinuity and uniformly in~$n$. The repaired law satisfies the
constraint exactly and is still close.

The proof of \cref{thm:existence} follows the strategy of Backhausz and Szegedy: from
the sequence, they construct a countable family of functions, closed under the
operators and under products, whose joint laws converge along a subsequence, and build
the limit map on a space whose coordinates are indexed by those functions.
Since nonlinear maps do not commute with linear combinations, we construct a larger
family, closed also under rational linear combinations, and hence under polynomials
with rational coefficients.


\section{Coupled maps on large graphs: convergence}\label{sec:applications}

In this section we apply \cref{thm:pairwise,thm:trajectories} to the coupled
maps~\eqref{eq:intro_nonlin} on large finite graphs, explain what the resulting
convergence means for the dynamics, and compute explicit limits for three examples,
from sparse to dense. Consequences for the dynamics are derived in \cref{sec:dynamics}.

\subsection{Coupled maps and their limits}\label{subsec:graph_dynamics}
We consider the coupled maps~$T_n$ in~\eqref{eq:network}, with a real matrix~$W_n$ and
continuous~$F,G\colon[-1,1]^2\to[-1,1]$ independent of~$n$. The function~$G$ describes
the influence of one vertex on another, while~$F$ combines the current state with the
aggregate interaction.

For each~$n$, we fix a probability measure~$\mu_n$ on~$[n]=\{1,\ldots,n\}$
with~$\mu_n(i)>0$ for every~$i$. The weight~$\mu_n(i)$ is the importance given to
vertex~$i$. There are some canonical choices for~$\mu_n$, and different choices can
describe different aspects of the dynamics, as we will see.
Assume that~$([n],\mu_n,W_n)\in\CL$, and that~$(W_n)$ is action convergent
and~$L^1$-equicontinuous. Explicitly, equicontinuity means that for every~$\ve>0$
there is~$\delta>0$, independent of~$n$, such that for all~$n$ and all~$u,v\in[-1,1]^n$,
\[
	\sum_{j=1}^n\mu_n(j)\,\abs{u(j)-v(j)}\leq\delta
	\quad\Longrightarrow\quad
	\sum_{i=1}^n\mu_n(i)\,\Bigl|\sum_{j=1}^nW_n(i,j)\bigl(u(j)-v(j)\bigr)\Bigr|\leq\ve.
\]
Requiring a uniform bound on~$\norm{W_n}_{1\to1}$ instead would be a strictly stronger
assumption, which would exclude interesting examples such as the scale-free trees
of \cref{subsec:scale_free}.

By \cref{thm:pairwise}, the sequence~$W_n$ has a limit~$W$, represented by
measures~$\nu_x$, and the sequence~$T_n=\Gamma(W_n)$ converges
to~$T:=\Gamma(W)$. Writing~$u(x,t)=(T^t(u))(x)$ for
the state of~$x\in\Omega$ at time~$t$, the limit of~$T_n$ reads
\begin{equation}\label{eq:limit_system}
	u(x,t+1)=F\Bigl(u(x,t),\ \int_\Omega G\bigl(u(x,t),u(y,t)\bigr)\d\nu_x(y)\Bigr),
	\qquad x\in\Omega.
\end{equation}

\subsubsection*{What convergence means}
We can now make precise the idea of \cref{subsec:intro_idea}.
Fix~$\ve>0$. For every~$n$ sufficiently large, the following is true:
for \textit{every}~$u\in[-1,1]^n$
there is~$w\in\Linf{\mu}$
such that~$d_W(\CD_{\mu_n}(u,T_nu),\CD_{\mu}(w,Tw))\leq \ve$,
and conversely for \textit{every}~$w\in\Linf{\mu}$
there is~$u\in[-1,1]^n$ such that~$d_W(\CD_{\mu_n}(u,T_nu),\CD_{\mu}(w,Tw))\leq \ve$.
This two-way approximation (which is the Hausdorff distance between profiles),
says that the limit captures the behavior of large finite systems and adds nothing
spurious to it.
By \cref{thm:trajectories}, the same holds for trajectories:
for every fixed~$t\geq 0$ and sufficiently large~$n=n(t,\ve)$, any law
\[
	\CD_{\mu_n}\bigl(u(\cdot,0),u(\cdot,1),\ldots,u(\cdot,t)\bigr)
	= \sum_{i=1}^n \mu_n(i)\, \delta_{(u(i,0),\ldots,u(i,t))}
\]
is $\ve$-close to some law
$
	\CD_\mu(w(\cdot,0),\ldots,w(\cdot,t)) = (x\mapsto (w(x,0),\ldots,w(x,t)))_\ast \mu
$ of the limit, and conversely.
More generally, same holds jointly for~$k$ initial conditions, for every~$k\geq 1$.

These laws carry dynamical information. If~$u$ is a fixed point of~$T_n$,
the law~$\CD_{\mu_n}(u,T_nu)$ is
supported on the diagonal.
If~$u$ is synchronized, that is~$u(1)=\cdots=u(n)$, the first
marginal is a Dirac mass; if the synchronized states also form an invariant set,
the law itself is a Dirac mass.
Note that~$\CD_{\mu_n}(u)=\CD_{\mu_n}(T_nu)$ does not
imply~$T_n u = u$: for instance, under uniform measure on the vertices,
a state and its permutations have the same law.
This is way considering the joint law~$\CD_{\mu_n}(u,T_nu)$ is important.

\subsection{Two normalizations}\label{subsec:normalizations}
There are several natural ways of representing a graph by a linear operator~$W_n$,
and several natural ways of choosing a probability measure~$\mu_n$ on its vertices.
These choices matter:
they determine which graph sequences are convergent, and the same graph sequence
may converge under different choices to different limits.
We discuss here two canonical choices. More are discussed in~\cite{backhausz2022action}.

\begin{example}[Uniform normalization]\label{ex:dense}
Let~$A_n$ be the adjacency matrix of a finite simple graph, and let~$W_n=A_n/n$
and~$\mu_n(i)=1/n$. Action convergence of~$(W_n)$ is equivalent to the usual dense
graph convergence~\cite[Section~8]{backhausz2022action}. The limit is given by a
graphon, a symmetric measurable function~$W\colon[0,1]^2\to[0,1]$: the
measures~$\nu_x$ in \cref{thm:pairwise} have density~$W(x,\cdot)$ with respect
to the Lebesgue measure on~$[0,1]$.
\end{example}

Represented as in \cref{ex:dense}, paths, cycles, trees, and stars
converge to the zero operator,
which is not very insightful. They have non-trivial limit
if represented as follows.

\begin{example}[Averaging over neighbours]\label{ex:randomwalk}
Let~$A_n$ be the adjacency matrix of a finite graph with no isolated vertices,
and let~$W_n(i,j)=A_n(i,j)/\deg(i)$, so that
\[
	(W_nu)(i)=\E_{j\sim i}u(j):=\frac1{\deg(i)}\sum_{j\sim i}u(j)
\]
is the average state over the neighbours of~$i$. Let~$\mu_n(i)=\deg(i)/\sum_j \deg(j)$ be
the degree measure. Action limits are then Markov graphops, and the
measures~$\nu_x$ in \cref{thm:pairwise} are probability
measures~\cite[Theorem~6.3]{backhausz2022action}.
\end{example}

The coupled maps~\eqref{eq:network} corresponding to \cref{ex:randomwalk} are
\begin{equation}\label{eq:examples_dynamics}
	(T_nu)(i)=F\Bigl(u(i),\ \E_{j\sim i}G\bigl(u(i),u(j)\bigr)\Bigr).
\end{equation}

We now present three graph sequences, from sparse to dense,
that are action convergent when represented as in \cref{ex:randomwalk}: cycles (bounded
degree), recursive scale-free trees (unbounded degree) and quasirandom graphs (dense).
We compute their limits explicitly. We also show how
changing~$\mu_n$ can change the limit: by giving more or less weight to different
parts of the graph, the limit captures different aspects of the dynamics.

\subsection{Example I: cycles}\label{subsec:cycles}
Let~$C_n$ be the cycle with vertices~$\Z/n\Z$ and edges~$\{i,i+1\}$, with~$n\ge3$,
so that~$(W_nu)(i)=\tfrac12\p{u(i-1)+u(i+1)}$, and let~$\mu_n(i)=1/n$.
Choose~$\alpha\in\R\setminus\Q$, and let~$\lambda$ be the Lebesgue measure
on~$\T:=\R/\Z$. Then~$(C_n,\mu_n,W_n)$ is action convergent to~$(\T,\lambda,W)$,
where~$(Wu)(x)=\tfrac12\p{u(x+\alpha)+u(x-\alpha)}$.
Indeed, this is the local--global limit of the cycles~\cite[Section~12]{hatami2014limits},
and for graphs of bounded degree local--global convergence is equivalent to action
convergence~\cite[Theorem~9.2]{backhausz2022action}.
Every irrational~$\alpha$ represents the same limit.
By \cref{thm:pairwise}, the corresponding
nonlinear maps~\eqref{eq:examples_dynamics}, namely
\[
	(T_nu)(i)=F\p{u(i),\ \tfrac12 (G(u(i),u(i-1)) + G(u(i),u(i+1)))},
\]
converge to
\[
	(Tu)(x)=F\Bigl(u(x),\tfrac12\p{G(u(x),u(x+\alpha))+G(u(x),u(x-\alpha))}\Bigr),
	\qquad x\in \T.
\]

The limit depends on the measure: changing~$\mu_n$ changes the importance of the
vertices, and hence the part of the dynamics that the limit describes.
For instance, fix~$p\in(0,1)$ and
let~$\mu_n(i)\propto p^{\abs{i}}$, where~$\abs{i}$ is the distance from~$i$ to~$0$
in~$C_n$. This essentially marks the vertex~$0$ and weights the other vertices by
their distance from it. Then~$(C_n,\mu_n,W_n)$ converges to~$(\Z,\mu,W)$,
where~$\mu(i)\propto p^{\abs{i}}$ and~$(Wu)(i)=\tfrac12\p{u(i-1)+u(i+1)}$.
Indeed, the operators
coincide for~$\abs{i}<n/2-1$, and the measures differ by~$\CO(p^{n/2})$ in total
variation. The nonlinear maps converge to
\[
 (Tu)(i)=F\Bigl(u(i),\tfrac12\p{G(u(i),u(i-1))+G(u(i),u(i+1))}\Bigr),
 \quad \ \ \qquad i\in \Z.
\]
With the uniform measure all vertices are indistinguishable, and the limit describes the
statistics of the whole cycle; with the geometric measure, the dynamics is seen from the ``point of view" of the
marked vertex~$0$. This also explains the two limit spaces:
$\T$ supports a translational-invariant probability measure
(making vertices indistinguishable),
$\Z$ doesn't.

With uniform measure~$\mu_n(i)=1/n$, the limit does not see changes on a vanishing fraction of
the vertices. For instance, removing the edge~$\{n,1\}$ turns~$C_n$ into the
path~$P_n$ and changes the operator only at two vertices, of total mass~$2/n$; similarly,
$C_n\sqcup C_n$ becomes~$C_{2n}$ after moving two edges. Hence,
with uniform measure, the sequences~$C_n$,
$P_n$, $C_{2n}$, $C_n\sqcup C_n$ all share the same action limit,
and so do their nonlinear dynamics.


\subsection{Example II: recursive scale-free trees}\label{subsec:scale_free}
The Barabási--Albert preferential attachment model is a standard model
of scale-free networks~\cite{barabasi1999emergence,bollobas2001degree}.
We consider a deterministic model for which a short argument gives
an explicit action limit, but our arguments also apply to the
Barabási--Albert model.

Start with an edge~$R_0$. To obtain~$R_{n+1}$, attach~$\deg_{R_n}(i)$ new leaves to
each vertex~$i$.
Thus, at each step,
every vertex doubles its degree.
This is the recursive scale-free tree of Jung,
Kim and Kahng~\cite{jung2002geometric} with~$m=1$; see also~\cite{zhang2009influences}.
A vertex born at step~$n-h$ has \emph{age}~$h$ and degree~$2^h$, and
\[
 |V(R_n)|=3^n+1,\qquad
 \#\{i:\operatorname{age}(i)=h\}=
 \begin{cases}2\cdot 3^{n-h-1},&0\le h<n,\\2,&h=n.\end{cases}
\]
We use the uniform measure~$\mu_n(i)=1/(3^n+1)$ and the operators~$W_n$ averaging
over neighbours in~$R_n$. We write~$\mathsf B_h$ for the tree
formed by a vertex of age~$h$ and its descendants.

To construct the limit, we label each vertex by its age~$h$ and a sequence~$x$ of
digits in~$\{0,1,2\}$, as follows. If the vertex has age~$h$ and its parent has
age~$h+j$, then the vertex is one of the~$2^{j-1}$ children that the parent received at
that step. We record which one by~$j-1$ digits in~$\{1,2\}$, followed by a~$0$, and
repeat this for the parent, the grandparent, and so on, up to~$R_0$.\footnote{For
example, a sequence starting with~$0020\ldots$ means that the vertex is the only child born
one step after its parent, the parent is likewise the only child born one step after the
grandparent, and the grandparent is child~$2$ of the two that the great-grandparent
received two steps after its own birth.}
The expected distance from~$R_0$ of a uniformly chosen vertex tends to infinity, so in the limit the sequences are infinite, as for the
cycle, where a uniformly chosen vertex is far from any fixed one. This motivates the
following.

Let~$\beta$ be the uniform product measure on~$\{0,1,2\}^{\N}$; sequences
with finitely many zeros form a null set, which we ignore.
Set~$\Omega=\N_0\times\{0,1,2\}^{\N}$ and~$\mu(\{h\}\times E)=\frac2{3^{h+1}}\,\beta(E)$.
The neighbours of a point~$(h,x)\in\Omega$ are defined as in~$R_n$. Write~$x=(w,0,y)$,
where~$w\in\{1,2\}^{j-1}$ consists of the digits before the first zero of~$x$ (empty
if~$x$ starts with~$0$). The \emph{parent} of~$(h,x)$ is~$(h+j,y)$, and its
\emph{children} are the points~$(h-k,(v,0,x))$ for~$1\le k\le h$
and~$v\in\{1,2\}^{k-1}$. There are one parent and~$\sum_{k=1}^{h}2^{k-1}=2^h-1$ children,
so~$(h,x)$ has~$2^h$ neighbours, and the operator that averages over neighbours is
\begin{equation}\label{eq:scale_uniform_operator}
 (Wu)(h,x)
 =2^{-h}\p{u(h+j,y)+\sum_{k=1}^{h}\ \sum_{v\in\{1,2\}^{k-1}}u\bigl(h-k,(v,0,x)\bigr)}.
\end{equation}
Note that this graph is scale-free: $\mu\{\deg\ge2^h\}=3^{-h}$.
One can show that, as~$n\to\infty$, the sequence~$(R_n,\mu_n,W_n)$ is action
convergent to~$(\Omega,\mu,W)\in\CT$; we omit the proof. Note that~$\norm{W_n}_{1\to1}\ge2^{n-1}$:
an initial vertex receives~$2^{n-1}$ leaves at the last step, and~$W_n$ maps its
indicator to~$1$ on these leaves. By \cref{thm:pairwise}, the
nonlinear maps~\eqref{eq:examples_dynamics} converge, and the limit is the following.
As before, for~$(h,x)\in\Omega$ write~$x=(w,0,y)$, where~$w\in\{1,2\}^{j-1}$ contains no
zeros and~$y\in\{0,1,2\}^{\N}$. Then
\[
 \begin{split}
 (Tu)(h,x)=F\Bigl(u(h,x),\ 2^{-h}\Bigl(&G\bigl(u(h,x),u(h+j,y)\bigr)\\
 &+\sum_{k=1}^{h}\ \sum_{v\in\{1,2\}^{k-1}}G\bigl(u(h,x),u(h-k,(v,0,x))\bigr)\Bigr)\Bigr).
 \end{split}
\]

As for cycles, changing the measure changes the part of the dynamics that the limit
describes. Here, we give positive measure to a vertex~$o$, a \textit{hub}.
In contrast to cycles,
the key difference is that here, as~$n\to\infty$,
the degree of~$o$ diverges.
Fix~$p\in(0,1)$, and let~$\tilde\mu_n$ give~$o$ mass proportional to~$p$, and the vertices at distance~$\ell$
from~$o$ total mass proportional to~$p^\ell$, split uniformly first among the
\emph{branches} (the connected components of~$R_n\setminus\{o\}$) that reach distance~$\ell$, and
then among the vertices of each such branch at distance~$\ell$.
For~$h<n$, the hub~$o$ has~$2^{n-h-1}$ branches rooted at its children of age~$h$, each
a copy of~$\mathsf B_h$, with root~$r_h$; as~$n\to\infty$ these copies form a continuum
labelled by~$x\in[0,1)$. Doing so gives a different limit, whose
points are the hub~$o$ and the points~$(h,z,x)$, where~$z$ is a vertex of~$\mathsf B_h$;
we omit the proof.
The nonlinear maps~\eqref{eq:examples_dynamics} converge to
\[
 \begin{aligned}
 (\tilde Tu)(o)&=F\Bigl(u(o),\ \sum_{h\ge0}2^{-h-1}\int_0^1G\bigl(u(o),u(h,r_h,x)\bigr)\d x\Bigr),\\
 (\tilde Tu)(h,z,x)&=F\Bigl(u(h,z,x),\ \frac1{\deg(z)}\Bigl(\sum_{z'\sim z}
 G\bigl(u(h,z,x),u(h,z',x)\bigr)\\
 &\hphantom{{}=F\Bigl(u(h,z,x),\ \frac1{\deg(z)}\Bigl(}
 +\1_{\{z=r_h\}}G\bigl(u(h,z,x),u(o)\bigr)\Bigr)\Bigr),
 \end{aligned}
\]
where~$z'\sim z$ denotes adjacency in~$\mathsf B_h$ and~$\deg(z)=2^{\operatorname{age}(z)}$
is the degree in~$R_n$. Compared with~$T$, the hub averages over a continuum of branch
roots, and each root also feels the hub. The limits of the cycle and of the uniform
measure are graphings; this one is not, since the hub has a continuum of neighbours.

\subsection{Example III: quasirandom graphs}\label{subsec:quasirandom}
\begin{figure}[ht]
  \centering
 \includegraphics[width=0.19\textwidth]{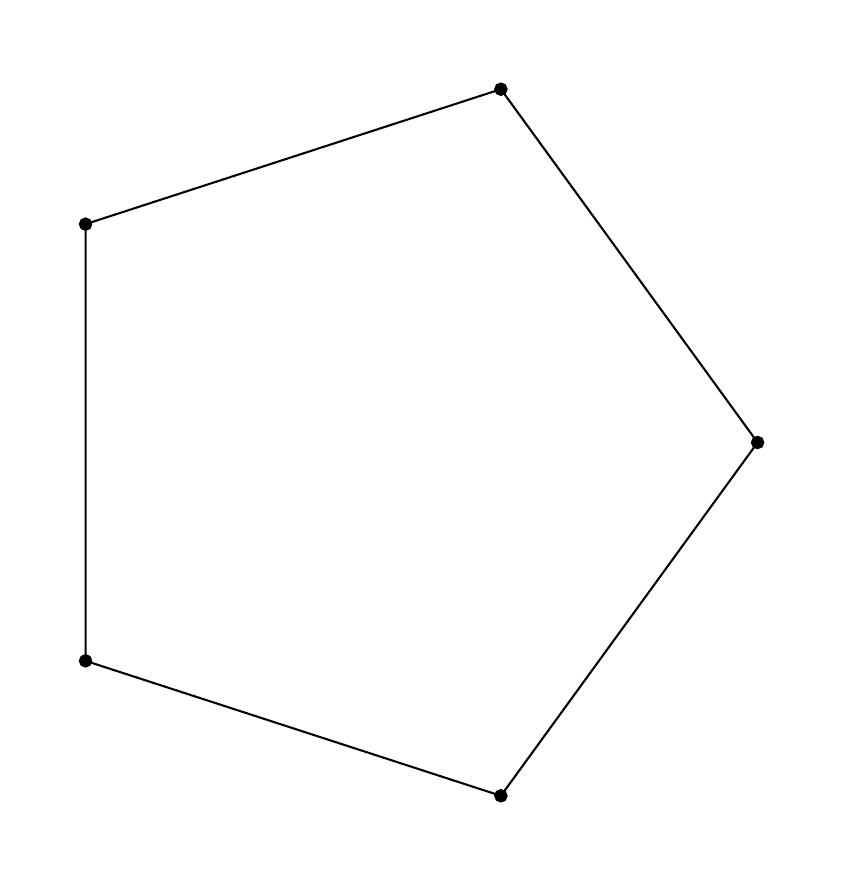}\hfill
 \includegraphics[width=0.19\textwidth]{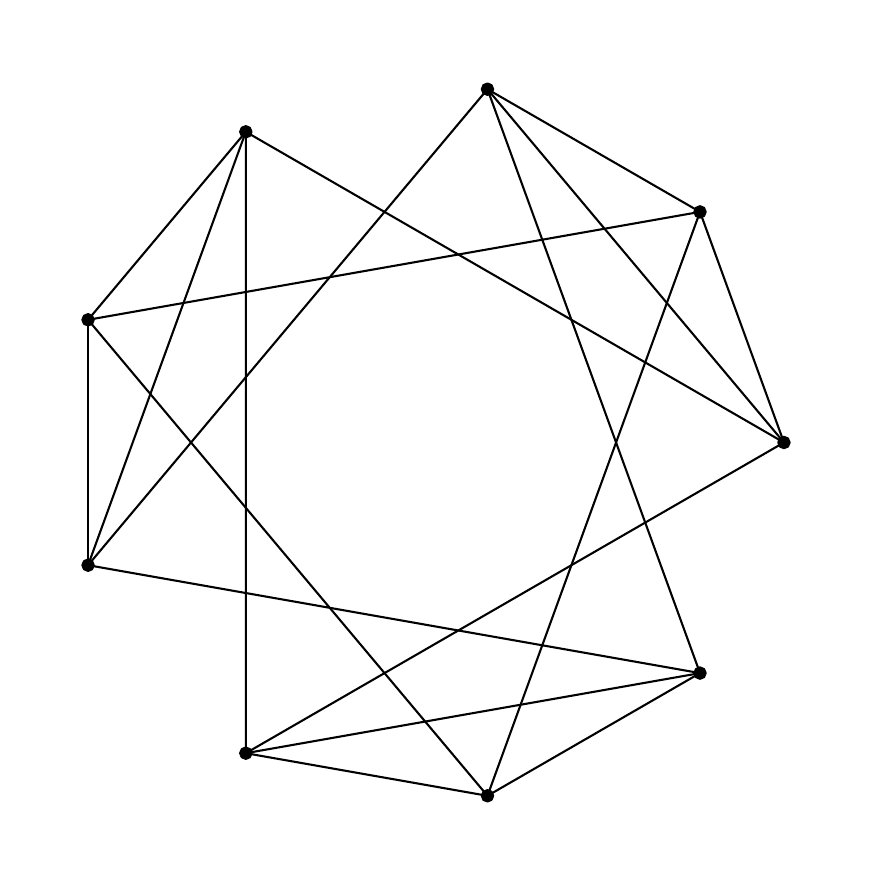}\hfill
 \includegraphics[width=0.19\textwidth]{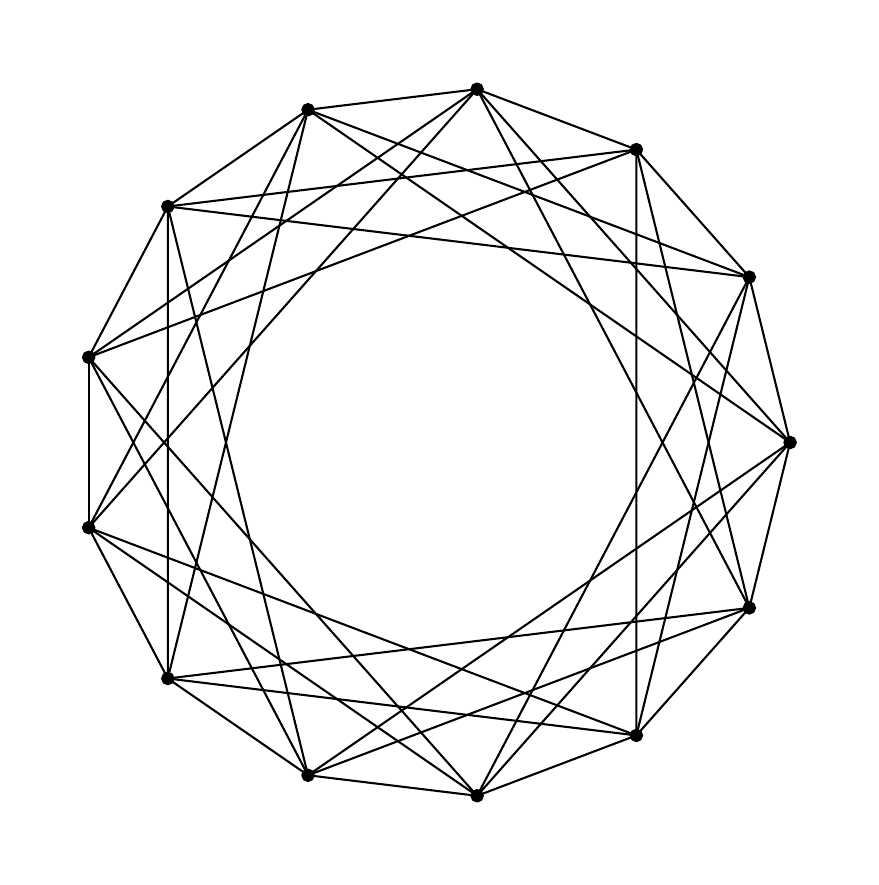}\hfill
 \includegraphics[width=0.19\textwidth]{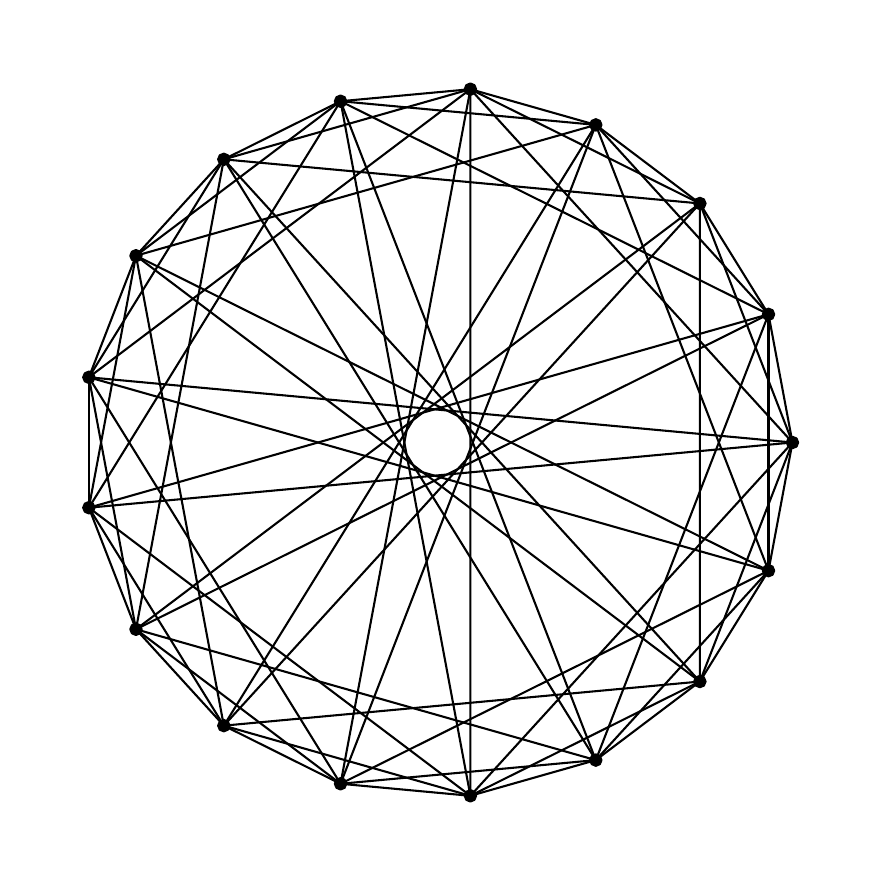}\hfill
 \includegraphics[width=0.19\textwidth]{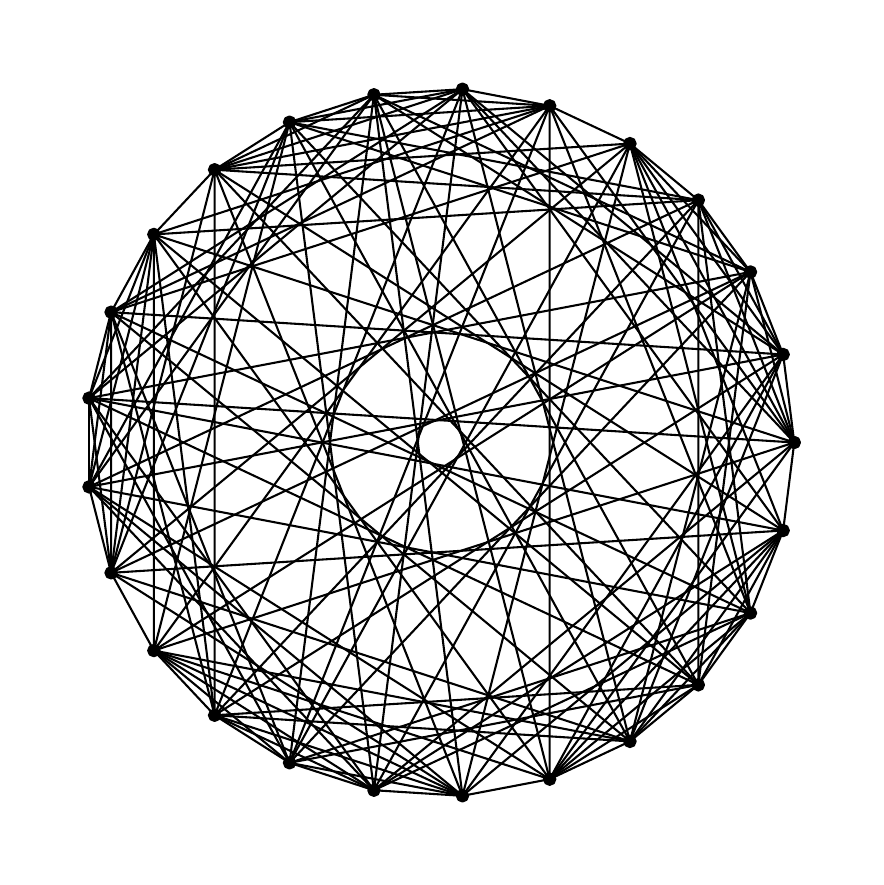}\\[1ex]
 \includegraphics[width=0.19\textwidth]{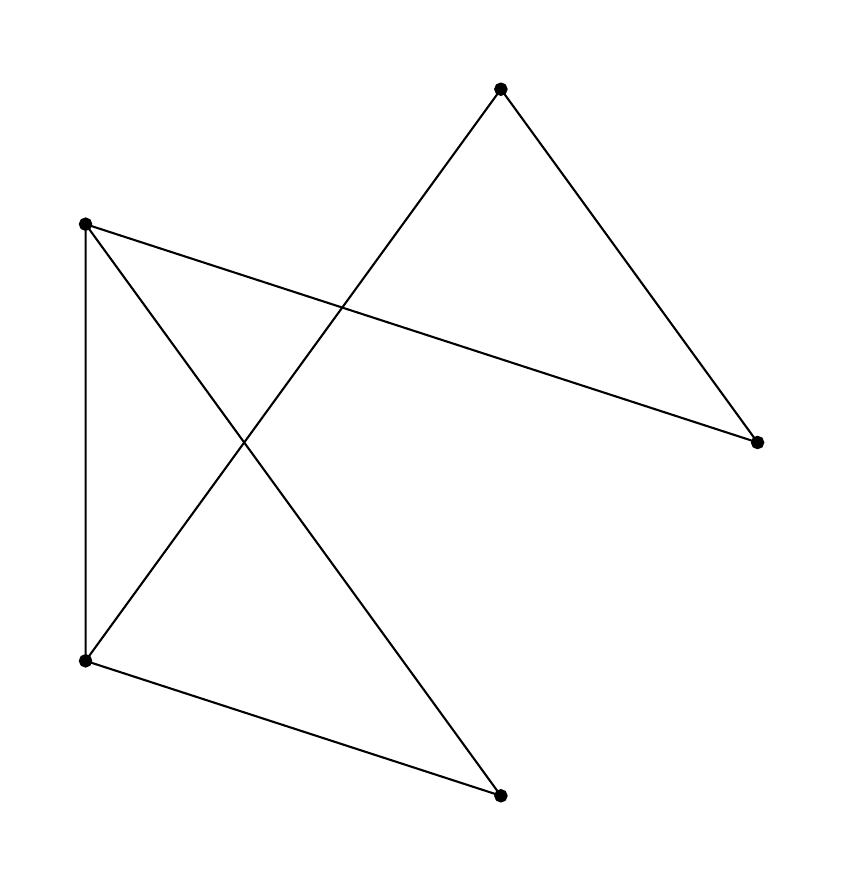}\hfill
 \includegraphics[width=0.19\textwidth]{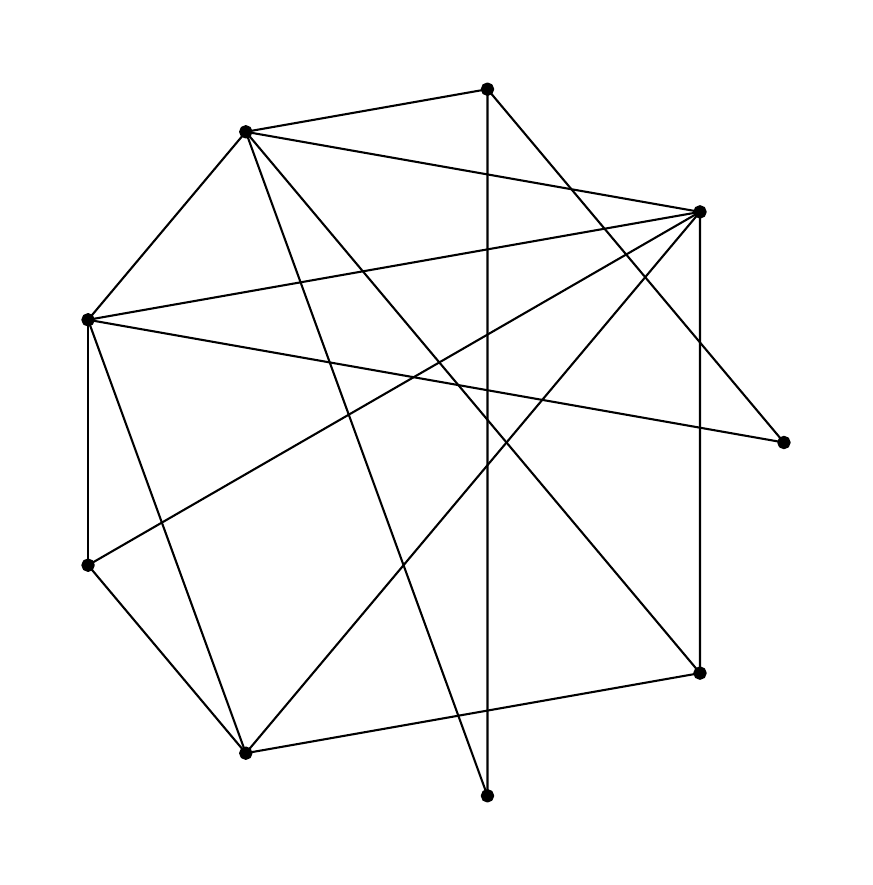}\hfill
 \includegraphics[width=0.19\textwidth]{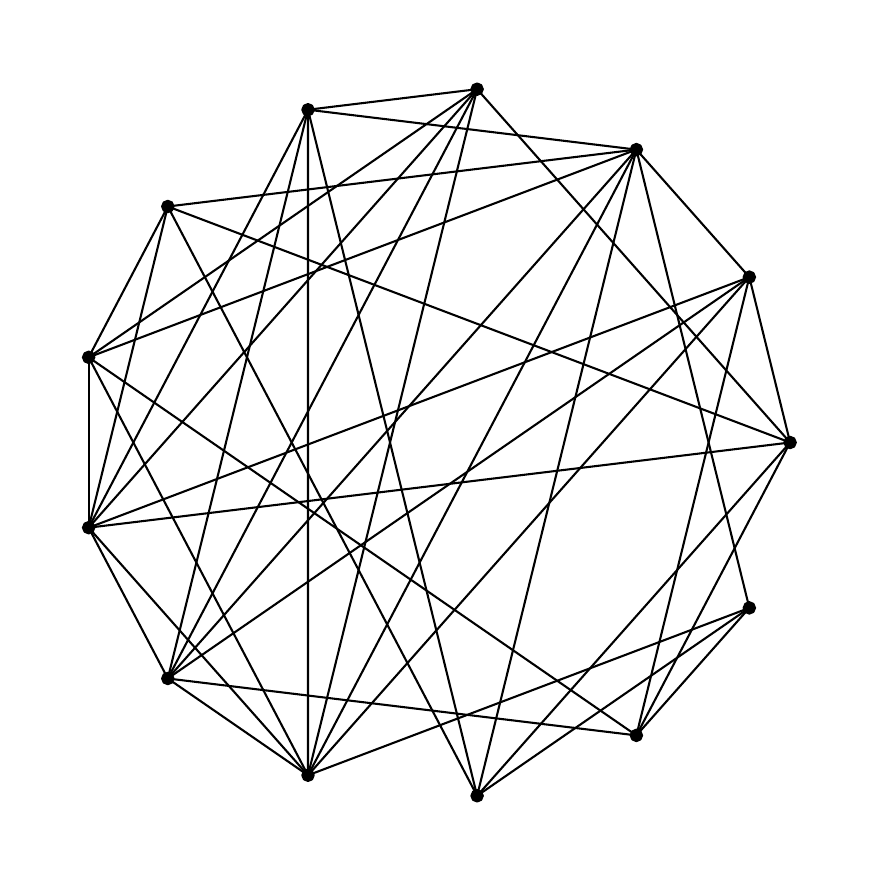}\hfill
 \includegraphics[width=0.19\textwidth]{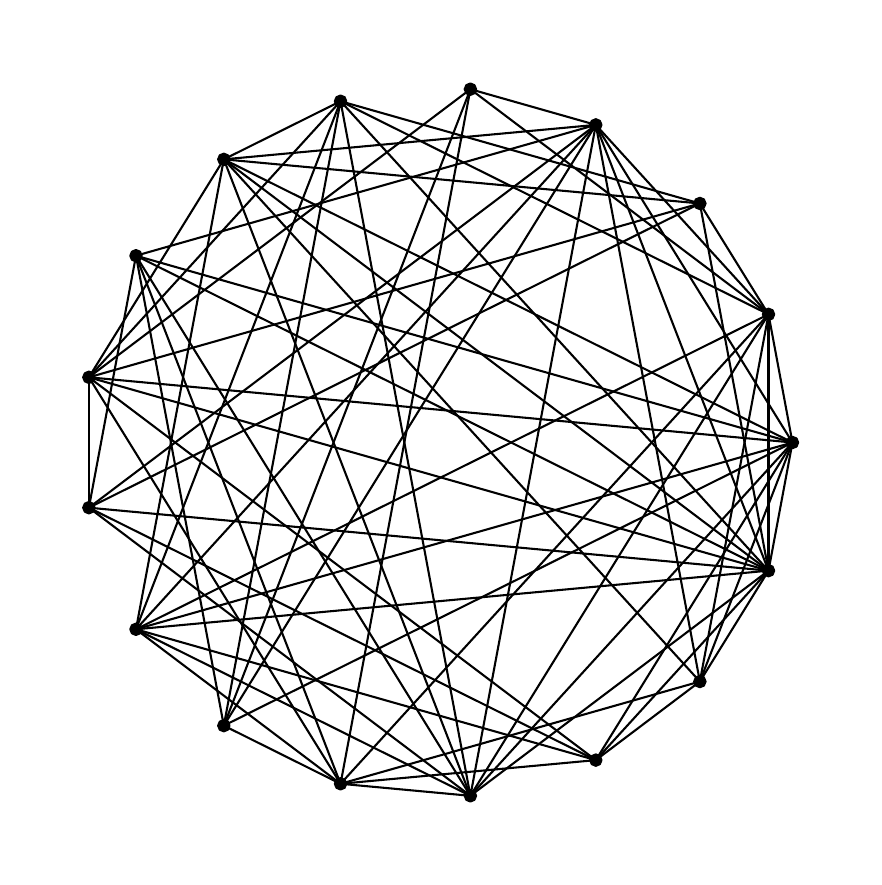}\hfill
 \includegraphics[width=0.19\textwidth]{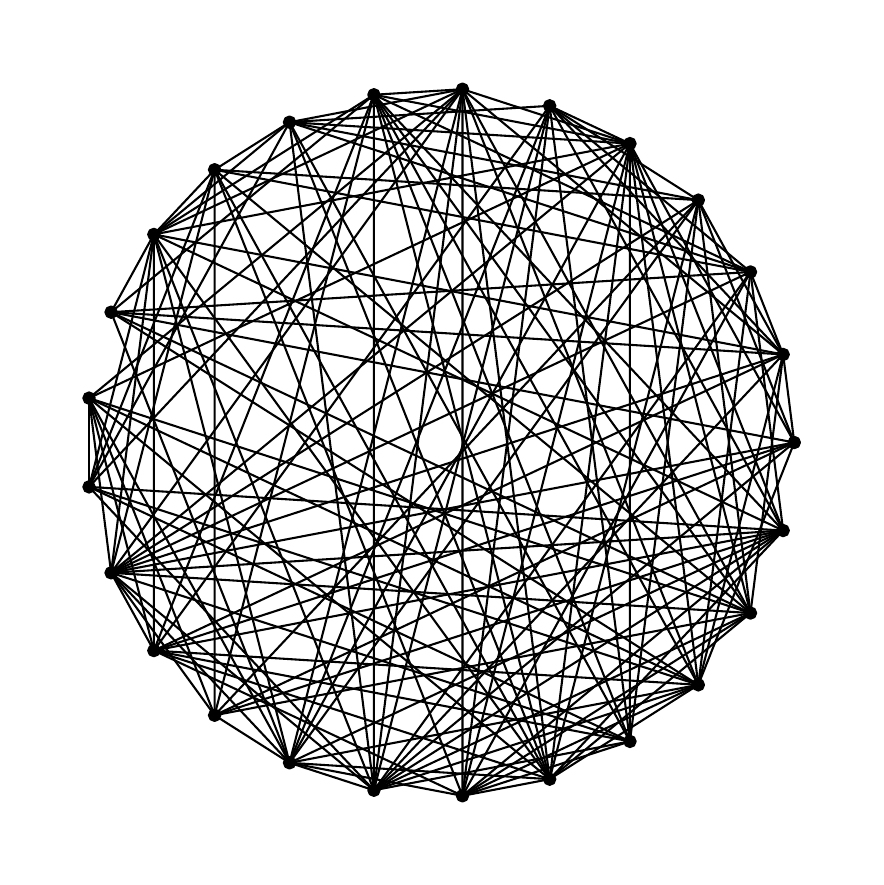}\\[1ex]
 \includegraphics[width=0.19\textwidth]{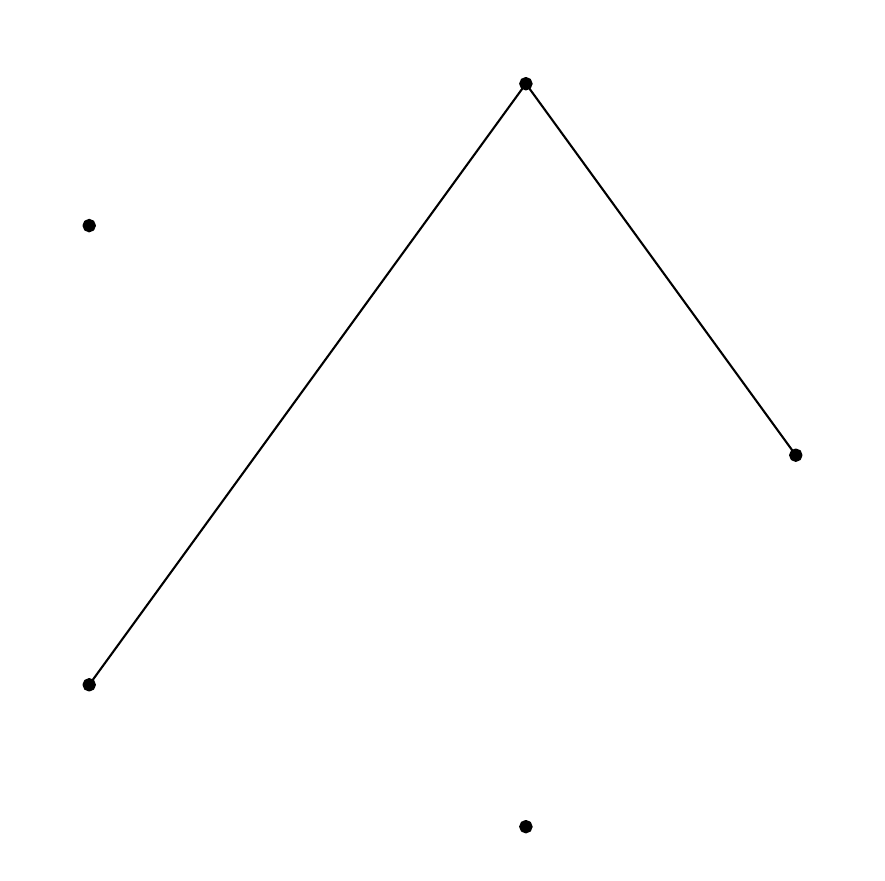}\hfill
 \includegraphics[width=0.19\textwidth]{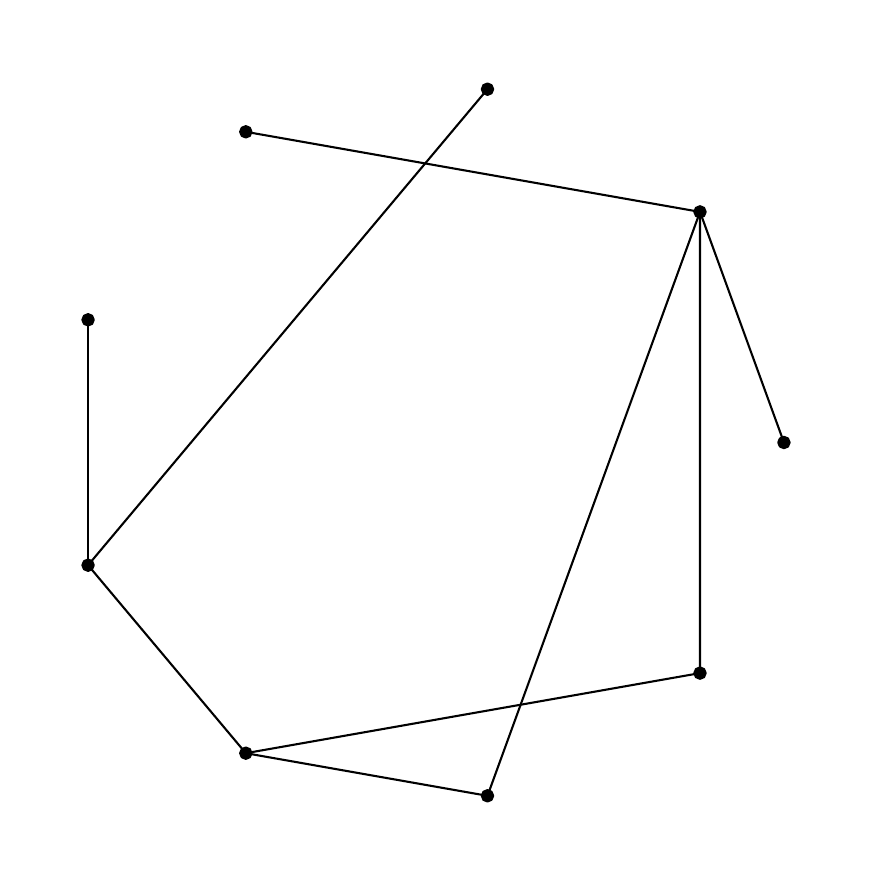}\hfill
 \includegraphics[width=0.19\textwidth]{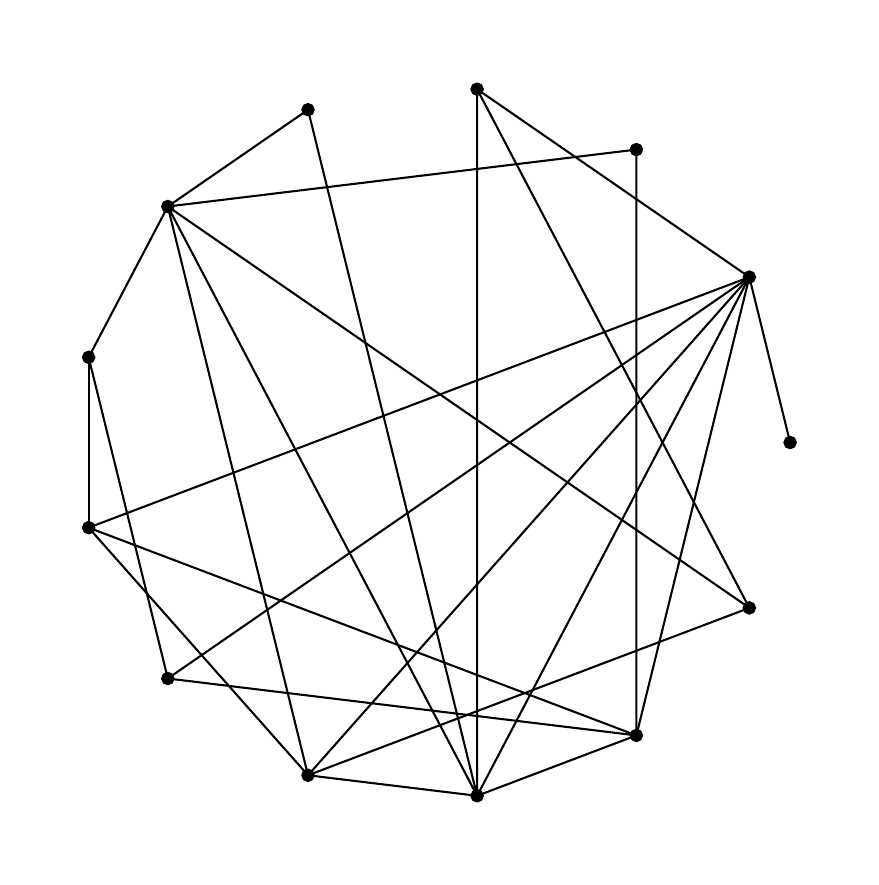}\hfill
 \includegraphics[width=0.19\textwidth]{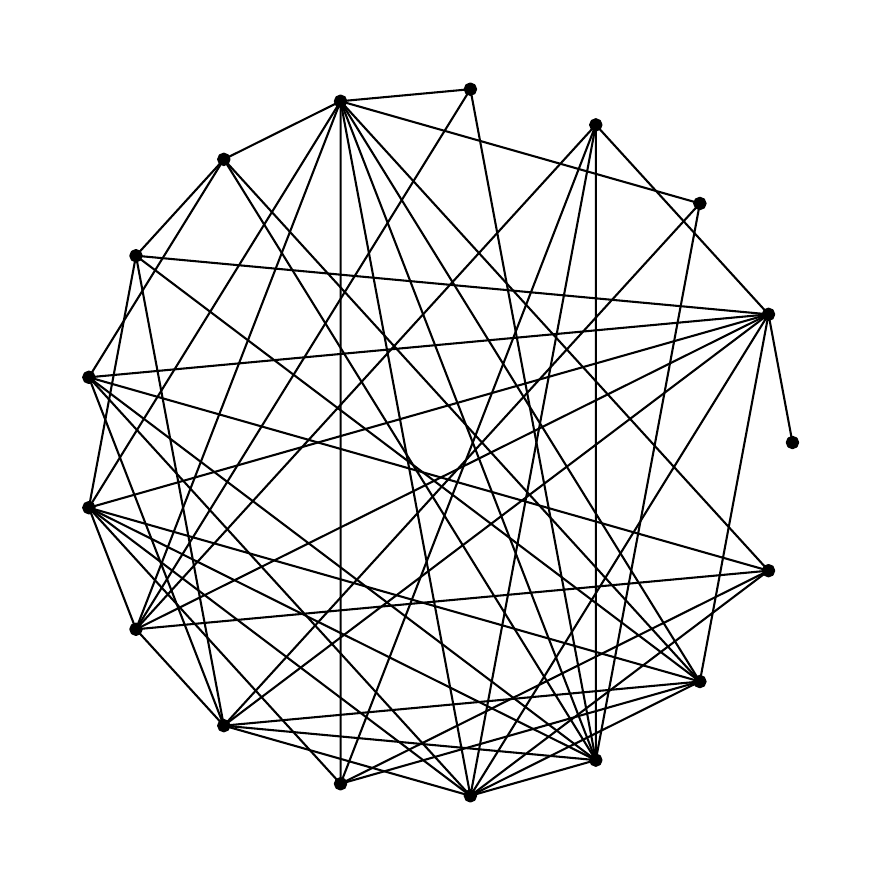}\hfill
 \includegraphics[width=0.19\textwidth]{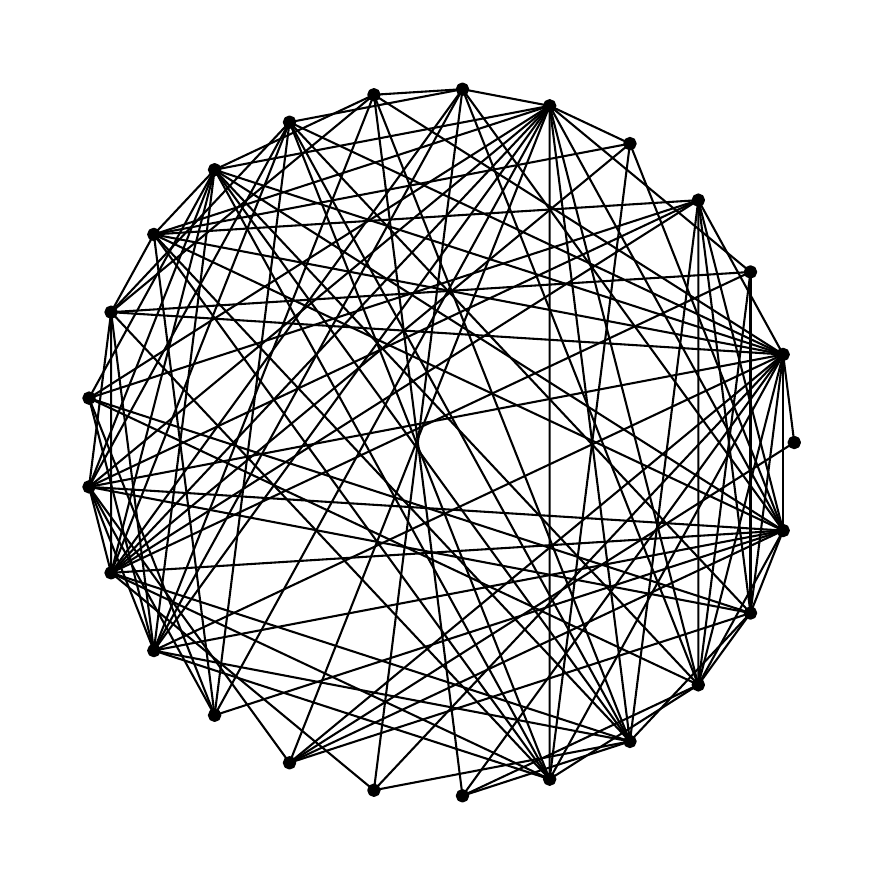}
  \caption{
  Top row: Paley graphs on $n = 5, 9, 13, 17, 25$ vertices.
	Middle row: Erd\H{o}s--R\'enyi random graphs $G(n,1/2)$ of the same orders.
	Bottom row: a single realization of the growing preferential attachment sequence,
	shown at $n = 5, 9, 13, 17, 25$ vertices.
  }
  \label{fig:er-paley}
\end{figure}

The previous examples are sparse: cycles have bounded degree, and the scale-free trees
have unbounded degrees but bounded average. At the other extreme are dense
graphs of which an important class are quasirandom graphs. A sequence of graphs~$(G_n)_{n\in\N}$, with~$G_n$ on~$[n]$ and possibly with loops, is
\emph{quasirandom} if there exists~$p\in(0,1)$ such that
\begin{equation} \label{eq:quasirandom}
	\Bigl\|\tfrac1n\mathrm{Adj}(G_n)- \frac{p}{n} \1_{n\times n} \Bigr\|_{\infty\to1}\to0,
\end{equation}
where $\mathrm{Adj}(G_n)$ denotes the adjacency matrix of $G_n$ and the norm is taken from~$L^\infty(\lambda_n)$ to~$L^1(\lambda_n)$;
equivalently, $G_n$ converges to the constant graphon~$p$.\footnote{
Associating an adjacency matrix to a graph requires ordering the vertices,
and $(\infty\to 1)$-convergence depends on the ordering in general
\cite[Theorem 11.59]{lovasz2012large}. This is irrelevant
for~\eqref{eq:quasirandom} due to
the full permutation symmetry of the matrix~$\frac{p}{n} \1_{n\times n}$.}
By the theorem of Chung,
Graham and Wilson~\cite{chung1989quasi}, this is equivalent to many other
conditions, for instance that almost all degrees are~$(p+o(1))n$ and almost all
codegrees, the numbers of common neighbours of two vertices, are~$(p^2+o(1))n$.
Examples of quasirandom graphs include:
\begin{itemize}
\item \emph{Paley graphs} (\cref{fig:er-paley}), with~$p=\frac12$~\cite{chung1989quasi};
\item \emph{strongly regular graphs} with parameters~$(n,d,a,b)$,
	where~$d\sim pn$ and~$a,b\sim p^2n$;
\item \emph{Erd\H{o}s--R\'enyi random graphs}~$G(n,p)$ with~$p$ fixed
	(\cref{fig:er-paley}): with probability~$1$, the sampled sequence is quasirandom;
\item \emph{growing preferential attachment graphs}~\cite[Example~11.44]{lovasz2012large}
	(\cref{fig:er-paley}). Start from a single vertex. At step~$n$, add a
	vertex~$n+1$ and join it to each~$i\in\{1,\ldots,n\}$ independently with
	probability~$(\deg_n(i)+1)/(n+1)$. With probability~$1$, the sampled
	sequence is quasirandom, for a value of~$p$ that depends on the realization;
\item \emph{perturbations} of any of these: changing~$o(n^2)$ edges, or adding or
	deleting~$o(n)$ vertices, preserves quasirandomness.
\end{itemize}

We assume that~$G_n$ has no isolated vertices, and use the operators~$W_n$ averaging over
neighbours and the degree measure~$\mu_n(i)=\deg(i)/\sum_j\deg(j)$, as in \cref{ex:randomwalk}. Since almost all
degrees are~$(p+o(1))n$, this measure is close to the uniform one. Let~$\lambda$ be the
Lebesgue measure on~$[0,1]$, and let~$J$ be the averaging operator,
$(Ju)(x)=\int_0^1u(y)\d y$, which replaces~$u$ by the constant function equal to its
mean. As~$n\to\infty$, the sequence~$([n],\mu_n,W_n)$ is action convergent
to~$([0,1],\lambda,J)\in\CT$.\footnote{By~\eqref{eq:quasirandom}
and~\cite[Section~8]{backhausz2022action}, the operators~$\frac1{pn}\mathrm{Adj}(G_n)$,
with the uniform measure~$\lambda_n$, converge to~$J$. Since almost all degrees
are~$(p+o(1))n$, one can verify by two approximations that replacing~$\frac1{pn}$
by~$\frac1{\deg(i)}$, and~$\lambda_n$ by the degree measure~$\mu_n$, does not change the
limit. We omit the details.} By \cref{thm:pairwise},
the nonlinear maps~\eqref{eq:examples_dynamics} converge to
\begin{equation}\label{eq:quasirandom_limit}
	(Tu)(x)=F\Bigl(u(x),\ \int_0^1G\bigl(u(x),u(y)\bigr)\d y\Bigr),
	\qquad x\in[0,1].
\end{equation}
Each vertex interacts with the whole system. Remarkably,
the limit depends neither on the
particular quasirandom sequence nor on~$p$. In \cref{subsec:STO} we show that this
limit is described exactly by the self-consistent transfer operator.


\section{Coupled maps on large graphs: dynamics}\label{sec:dynamics}

In \cref{sec:applications} we described the limit objects of coupled maps on large
graphs. Here we show how to use them to understand the dynamics of large but finite
networks. The dynamical properties we obtain are, by construction, robust under small
changes of the graph, which do not affect the limit, and hold uniformly over all initial
conditions. Most of them concern finite time horizons, but some are uniform in time.

\subsection{The self-consistent transfer operator and quasirandom graphs}
\label{subsec:STO}
Fix continuous functions~$F,G\colon[-1,1]^2\to[-1,1]$.
The \emph{self-consistent transfer operator} (STO) is the continuous map
$S\colon\CM([-1,1])\to\CM([-1,1])$ defined by
\begin{equation}\label{eq:STO}
	S(\nu):=(\Phi_\nu)_*\nu,
	\qquad
	\Phi_\nu(z):=F\Bigl(z,\ \int G(z,y)\d\nu(y)\Bigr).
\end{equation}
For a fixed measure~$\eta$, the push-forward~$\nu\mapsto(\Phi_\eta)_*\nu$ is the
linear transfer operator of~$\Phi_\eta$. Letting the map depend on the measure it
acts on makes~$S$ nonlinear in general. In the
literature~\cite{keller2000ergodic,balint2018synchronization,selley2021linear,galatolo2022self},
the maps~$\Phi_\nu$ are usually smooth and nonsingular, so that~$S$ preserves
absolutely continuous measures and can be studied as an operator on densities.
Here we require only continuity of~$F$ and~$G$ and work with all probability
measures, including atomic empirical distributions.

It is well known that the STO describes the evolution of empirical distributions
in all-to-all coupled systems with full permutation symmetry, i.e.\ in
equation~\eqref{eq:network} with~$W_n(i,j)=1/n$ for all~$i,j$~\cite{kaneko1990globally,keller2000ergodic}.
Indeed, for a state~$u\in[-1,1]^n$, the empirical distribution is precisely its
law under the uniform vertex measure:
$
	\nu_u:=\frac1n\sum_{i=1}^n\delta_{u(i)}=\CD_{\lambda_n}(u).
$
The update~\eqref{eq:network} then reads~$T_nu=\Phi_{\nu_u}\circ u$, hence
\[
	\CD_{\lambda_n}(T_n^tu)=S^t\bigl(\CD_{\lambda_n}(u)\bigr)
	\qquad\text{for every~$t\geq0$}.
\]
For every~$n$, the map~$u\mapsto \CD_{\lambda_n}(u)$
is a factor map from the coupled system to the STO.

On a general network, the empirical distribution of a state does
not determine its own evolution.
On the cycle~$C_n$
with~$4\mid n$, take~$W_nu(i)=\frac12\bigl(u(i-1)+u(i+1)\bigr)$, and consider
the two initial conditions
\[
	(+1,+1,-1,-1,+1,+1,-1,-1,\ldots)
	\qquad\text{and}\qquad
	(+1,\ldots,+1,-1,\ldots,-1).
\]
Both have distribution~$\frac12(\delta_{-1}+\delta_{+1})$. After one step, the first has
distribution~$\delta_0=S\bigl(\frac12(\delta_{-1}+\delta_{+1})\bigr)$, while the second is
unchanged except at four vertices. The distance between these distributions remains
bounded away from~$0$ as~$n\to \infty$. The STO does not describe the system.
We show that,
nevertheless, the STO is asymptotically exact on quasirandom graphs.
This may seem surprising,
as quasirandomness is a deterministic condition and does not require permutation
symmetry of the finite graphs. In fact, quasirandom graphs may have
no symmetries at all other than the identity.

Consider the dynamics~\eqref{eq:examples_dynamics},
with the degree measure~$\mu_n$ and the operator~$W_n$ averaging over
neighbours (\cref{ex:randomwalk}).
Let~$(G_n)$ be quasirandom, as defined in~\cref{subsec:quasirandom}:
the operators~$W_n$ converge to the averaging
operator~$J$, and the maps~$T_n$ to the limit system~\eqref{eq:quasirandom_limit}
on~$([0,1],\lambda)$. The observation is that this limit
factors through the STO:
\[
	\CD_\lambda(Tu)=S\p{\CD_\lambda(u)}\qquad\text{for all~$u\in\Linf{\lambda}$}.
\]
Action convergence transfers it to the finite graphs.

\begin{corollary}\label{cor:STO}
If~$(G_n)$ is quasirandom, then there is a sequence~$\ve_n\to0$ such that, for all~$n$
and all~$u\in[-1,1]^n$, the maps~$T_n$ in~\eqref{eq:examples_dynamics} satisfy
\[
	d_W\bigl(\CD_{\mu_n}(T_nu),\,S(\CD_{\mu_n}(u))\bigr)\leq \ve_n.
\]
\end{corollary}

\begin{proof}
Since~$T_n\to T$, we have~$\delta_n:=d_H\bigl(P_1(T_n),P_1(T)\bigr)+1/n\to0$.
Given~$u\in[-1,1]^n$, choose~$w\in\Linf{\lambda}$ with
$d_W\bigl(\CD_{\mu_n}(u,T_nu),\CD_\lambda(w,Tw)\bigr)\leq\delta_n$.
Taking marginals, which is $1$-Lipschitz for~$d_W$, and using~$\CD_\lambda(Tw)=S(\CD_\lambda(w))$,
we get~$d_W(\CD_{\mu_n}(u),\CD_\lambda(w))\leq\delta_n$
and~$d_W\bigl(\CD_{\mu_n}(T_nu),S(\CD_\lambda(w))\bigr)\leq\delta_n$.
Let~$\kappa$ be a modulus of uniform continuity of~$S$ on the compact
space~$\CM([-1,1])$. Then the claim holds with~$\ve_n:=\delta_n+\kappa(\delta_n)$.
\end{proof}


\subsection{Approximate fixed points}\label{subsec:fixed}
A state~$u$ is a \emph{fixed point} of~$T$ if~$Tu=u$. For~$\nu\in\CM([-1,1])$,
let~$\hat\nu$ be its push-forward under~$a\mapsto(a,a)$. Then~$u$ is a fixed point
with law~$\nu$ if and only if~$\CD_\mu(u,Tu)=\hat\nu$, and since
$\norm{Tu-u}_{L^1(\mu)}=\int\abs{a-b}\d\CD_\mu(u,Tu)$, a sequence of laws~$\CD_\mu(u_k,Tu_k)$
converges to~$\hat\nu$ if and only if~$\CD_\mu(u_k)\to\nu$
and~$\norm{Tu_k-u_k}_{L^1(\mu)}\to0$.

\begin{example}[Cycles, paths and disjoint cycles] \label{ex:cycles}
On the cycle~$C_n$ consider the dynamics
\[
	T_nu(i)=\tfrac12\bigl(G(u(i),u(i-1))+G(u(i),u(i+1))\bigr),
\]
that is, \eqref{eq:examples_dynamics} with~$F(a,c)=c$, and its analogues
on paths~$P_n$ and on disjoint unions of cycles~$C_n\sqcup C_n$. The graph~$C_{2n}$
becomes~$P_{2n}$ after removing one edge, and~$C_n\sqcup C_n$ after moving two edges,
while~$C_{2n+1}$ has just one more vertex. As observed in \cref{subsec:cycles}, all four
families have the same limit, and so do the
corresponding dynamics~$T_n$. Their fixed points, however, differ:
\begin{itemize}
\item if~$G(a,a)=a$ and~$G(b,b)=b$, the state equal to~$a$ on one cycle of~$C_n\sqcup C_n$
	and to~$b$ on the other is fixed; on~$C_{2n}$ or~$P_{2n}$, the corresponding
	two-arc state is in general not fixed, because the vertices at the interfaces
	see both values, but~$\norm{T_nu-u}_{L^1(\mu_n)}=\CO(1/n)$;
\item if~$G(a,b)=a$ and~$G(b,a)=b$, the alternating state~$(a,b,a,b,\ldots)$ is fixed
	on~$C_{2n}$; on~$C_{2n+1}$ the alternation breaks at one vertex,
	and~$\norm{T_nu-u}_{L^1(\mu_n)}=\CO(1/n)$.
\end{itemize}
So whether~$\frac12(\delta_a+\delta_b)$ is the law of a fixed point depends
on~$\CO(1)$ edges: it is not an asymptotic property of the system.
\end{example}

Exact fixed points therefore cannot be transferred to, or from, the limit.
Approximate fixed points can, in both directions.

\begin{corollary} \label{cor:fixed}
Let~$T_n\to T$ and~$\nu\in\CM([-1,1])$. The
following are equivalent:
\begin{enumerate}[(i)]
\item there are states~$w_k\in\Linf{\mu}$ with~$\CD_\mu(w_k)\to\nu$
	and~$\norm{Tw_k-w_k}_{L^1(\mu)}\to0$;
\item there are states~$u_n\in\Linf{\mu_n}$ with~$\CD_{\mu_n}(u_n)\to\nu$
	and~$\norm{T_nu_n-u_n}_{L^1(\mu_n)}\to0$.
\end{enumerate}
\end{corollary}

\begin{proof}
By the observation above, applied on each space, (i) means
that~$\hat\nu\in\overline{P_1(T)}$, and~(ii) means that there
are~$\rho_n\in P_1(T_n)$ with~$\rho_n\to\hat\nu$.
Since~$\overline{P_1(T_n)}\to \overline{P_1(T)}$,
the two conditions are equivalent.
\end{proof}


\subsection{Macroscopic Lyapunov functions, synchronization, and clusters}
\label{subsec:lyapunov}

To obtain bounds uniform in time, we look for Lyapunov functions of the limit that
depend only on the distribution of states, and are therefore seen by the profiles.

\begin{definition}
Fix~$\theta\in[0,1)$ and a function~$L\colon\CM([-1,1])\to[0,1]$, continuous in the
Wasserstein metric. We say that~$T$ is \emph{$L$-contracting with rate~$\theta$} if
\[
	L(\CD_\mu(Tu))\leq \theta\,L(\CD_\mu(u))
	\qquad\text{for all~$u\in\Linf{\mu}$,}
\]
or, equivalently,
if~$L(\pi_2\rho)\leq \theta\,L(\pi_1\rho)$
for all~$\rho\in\overline{P_1(T)}$,
where~$\pi_1,\pi_2$ denote the marginals.
\end{definition}

The second formulation shows that $L$-contraction depends only on the profile,
which makes it well suited to action convergence. By compactness
of~$\CM([-1,1])$, any continuous~$L\geq0$ can be normalized to take values
in~$[0,1]$.

\begin{theorem} \label{thm:Lyapunov}
Let~$T_n\to T$. Then~$T$ is $L$-contracting with rate~$\theta$ if and only if there
are~$\ve_n\geq0$ with~$\ve_n\to0$ such that
\begin{equation} \label{eq:Lyapunov1}
	L(\CD_{\mu_n}(T_nu))\leq \theta\,L(\CD_{\mu_n}(u))+\ve_n
	\qquad\text{for all~$u\in\Linf{\mu_n}$}.
\end{equation}
In that case, for all~$n$, all~$u\in\Linf{\mu_n}$ and all~$t\geq0$,
\begin{equation} \label{eq:Lyapunov2}
	L(\CD_{\mu_n}(T_n^tu))\leq \theta^t L(\CD_{\mu_n}(u))+\frac{\ve_n}{1-\theta}.
\end{equation}
\end{theorem}

\begin{proof}
For a compact~$K\subseteq\CM([-1,1]^2)$, let
$\tilde L(K):=\{(L(\pi_1\rho),L(\pi_2\rho))\colon\rho\in K\}\subseteq[0,1]^2$.
Since~$L\circ\pi_1$ and~$L\circ\pi_2$ are uniformly continuous, $\tilde L(K)$ is
compact and~$K\mapsto\tilde L(K)$ is continuous for the Hausdorff metric; hence
$\tilde L(\overline{P_1(T_n)})\to\tilde L(\overline{P_1(T)})$. For~$\ve\geq0$, let
$\Delta^\ve:=\{(x,y)\in[0,1]^2\colon y\le \theta x+\ve\}$, a closed set.
Then~$T$ is $L$-contracting with rate~$\theta$ if and only
if~$\tilde L(\overline{P_1(T)})\subseteq\Delta^0$, which by convergence is equivalent
to~$\tilde L(\overline{P_1(T_n)})\subseteq\Delta^{\ve_n}$ for some~$\ve_n\to0$,
that is, to~\eqref{eq:Lyapunov1}. For~\eqref{eq:Lyapunov2}, apply~\eqref{eq:Lyapunov1}
to~$T_n^su$ for~$s=0,\ldots,t-1$ and use~$1+\theta+\cdots+\theta^{t-1}\le(1-\theta)^{-1}$.
\end{proof}

The basic example is the variance~$\var(\nu):=\int x^2\d\nu(x)-\p{\int x\d\nu(x)}^2$.
It is continuous on~$\CM([-1,1])$, takes values in~$[0,1]$, and vanishes exactly at
the Dirac measures~$\delta_c$, the laws of synchronized states.

\begin{corollary}[Synchronization] \label{cor:synch}
Let~$T_n\to T$, and suppose that~$T$ is $\var$-contracting with rate~$\theta$.
Then there are~$\ve_n\to0$ such that, for all~$n$, all~$u\in\Linf{\mu_n}$ and
all~$t\geq0$,
\begin{equation} \label{eq:synch}
	d_W\bigl(\CD_{\mu_n}(T_n^tu),\delta_{c(t)}\bigr)\leq \theta^{t/2}+\ve_n,
\end{equation}
where~$c(t)$ is the $\mu_n$-mean of~$T_n^tu$.
\end{corollary}

\begin{proof}
By \cref{thm:Lyapunov}, there are~$\ve'_n\to0$ such
that~$\var(\CD_{\mu_n}(T_n^tu))\leq \theta^t+\ve'_n$ for all~$n$, $u$ and~$t$.
For~$\nu\in\CM([-1,1])$ with mean~$m$, Cauchy--Schwarz gives
$d_W(\nu,\delta_m)=\int\abs{x-m}\d\nu\leq\sqrt{\var(\nu)}$.
Hence~\eqref{eq:synch} holds with~$\ve_n:=\sqrt{\ve'_n}$,
since~$\sqrt{a+b}\leq\sqrt a+\sqrt b$.
\end{proof}

\begin{example}[Synchronization on quasirandom graphs] \label{ex:synch}
Let~$G_n$ be quasirandom, with the degree measures~$\mu_n$ as in \cref{subsec:quasirandom}.
Let~$F\colon[-1,1]\to[-1,1]$ be Lipschitz, let~$h\in[0,1]$, and consider
$
	(T_nu)(i):=(1-h)\,F(u(i))+h\,\E_{j\sim i}F(u(j)).
$
We have~$T_n\to T$, where~$Tu=(1-h)F\circ u+h\int F\circ u\d\lambda$. Since~$Tu$
is~$(1-h)F\circ u$ plus a constant,
$
	\var(\CD_\lambda(Tu))=(1-h)^2\var(\CD_\lambda(F\circ u))\leq\theta\,\var(\CD_\lambda(u)),
$
with~$\theta:=(1-h)^2\lip{F}^2$. If~$(1-h)\lip{F}<1$, then~$T$ is $\var$-contracting,
and \cref{cor:synch} gives
\begin{equation} \label{eq:synch_quasirandom}
	\sum_{i=1}^n\mu_n(i)\,\abs{u(i,t)-c(t)}\leq \theta^{t/2}+\ve_n
	\qquad\text{for all~$t\geq0$ and all~$u\in[-1,1]^n$},
\end{equation}
where~$c(t):=\sum_i\mu_n(i)u(i,t)$.
In words: on every large quasirandom graph, after a transient whose length is
independent of the initial condition, all but a vanishing fraction of the vertices are
close to the common value~$c(t)$, and they remain so at all later times.
On complete graphs, synchronization follows from a direct computation, since all
vertices see essentially the same average. Here, instead,~$G_n$ can be any of the
graphs of \cref{fig:er-paley}, or a perturbation of them, and a direct proof
of~\eqref{eq:synch_quasirandom} would need graph-specific estimates.
\end{example}

\begin{remark}[Clusters] \label{rem:clusters}
Since~$\var(\nu)=\min_a\int\abs{x-a}^2\d\nu(x)$, \cref{cor:synch} controls the
distance from laws supported on one point, i.e. full synchrony.
For two clusters, take
\[
	L(\nu):=\min_{a,b\in[-1,1]}\int\min\bigl\{\abs{x-a}^2,\abs{x-b}^2\bigr\}\d\nu(x).
\]
It is continuous, takes values in~$[0,1]$, and vanishes exactly at the laws supported on at most two points.
The same idea extends to any number of clusters.
\end{remark}

\begin{remark}[Microscopic Lyapunov functions]
The functions~$L$ are \emph{macroscopic}: they see only the law~$\CD_\mu(u)$. Together
with the additive error in~\eqref{eq:Lyapunov1}, this makes \cref{thm:Lyapunov}
insensitive to changes of the graph on a vanishing set of vertices. Exact contraction
is not robust in this sense. Let~$W_n$ be the linear
operator averaging over neighbours, with the degree
measure, on~$K_n\sqcup K_2$, where~$K_n$ is the complete graph. The component~$K_2$ has
vanishing mass and disappears in the limit, which is $\var$-contracting
with~$\theta=0$. Yet the state equal to~$0$ on~$K_n$ and to~$1,-1$ on~$K_2$ is
$2$-periodic but not fixed. This rules out any \emph{microscopic} Lyapunov function:
if~$L\colon [-1,1]^{n+2} \to [0,1]$
satisfies~$L(W_n u)\leq\theta_nL(u)$ for some~$\theta_n<1$
and vanishes exactly at the fixed points, then every
$2$-periodic point must be a fixed point.
For the same reason, \eqref{eq:Lyapunov1} fails with~$\ve_n=0$.
\end{remark}


\section{Proofs of \texorpdfstring{\cref{thm:pairwise,thm:trajectories}}{Theorems B and C}}\label{sec:lifting}

In this section we prove the ``lift and repair'' result announced in the
introduction (\cref{lem:repair,cor:lifting}) and deduce from it
\cref{thm:trajectories}, the case of a single inner function
(\cref{prop:single}), and \cref{thm:pairwise}. We end with a converse
of \cref{thm:pairwise}.

Let~$K\subseteq[-1,1]^{2m}$ be closed and~$\CM(K)\subseteq \CM([-1,1]^{2m})$
the set of measures supported on~$K$.
In general, $P_m(T_n)\to P_m(T)$
does not imply $P_m(T_n)\cap\CM(K)\to P_m(T)\cap\CM(K)$.
It does, however,
when~$K$ is defined by a triangular system of equations. Let us make this precise.
Points of~$[-1,1]^{2m}$ are written~$z=(x,y)$, with~$x=(x_1,\ldots,x_m)$
and~$y=(y_1,\ldots,y_m)$.
We write
$
	z_{<j}:=(x_1,\ldots,x_{j-1},y_1,\ldots,y_{j-1})\in[-1,1]^{2(j-1)}.
$
For a map~$T\colon\Linf{\mu}\to\Linf{\mu}$
and~$u\in\Linf{\mu}^m$, the coordinates~$x_j$ and~$y_j$ of~$\CD_\mu(u,Tu)$
record the input~$u_j$ and the output~$Tu_j$.

\begin{definition}[Triangular constraint] \label{def:triangular}
A \emph{triangular constraint}~$g$ on~$[-1,1]^{2m}$ consists of a
set~$J\subseteq\{1,\ldots,m\}$ and continuous maps~$g_j\colon[-1,1]^{2(j-1)}\to[-1,1]$
for~$j\in J$. Define
\[
	K_g:=\{z\in[-1,1]^{2m}\colon x_j=g_j(z_{<j})\text{ for all }j\in J\},
	\qquad
	\varphi_g(z):=\sum_{j\in J}\abs{x_j-g_j(z_{<j})}.
\]
\end{definition}

Each constrained input~$x_j$ is a continuous function of the inputs
and outputs of earlier coordinates, never of its own output or of later ones.
In the applications, the $k$-profile of a map is analyzed by constraining
the $m$-profile of a possibly different map, for some~$m \geq k$.

\begin{table}[ht]
\centering
\small
\caption{Triangular constraints used in this paper. The second column
describes the constrained profile~$P_m(T)\cap\CM(K_g)$:
the configurations~$u\in \Linf{\mu}^m$
for which~$\CD(u,Tu)$ is supported on~$K_g$.}
\label{tab:constraints}
\begin{tabular}{@{}llll@{}}
\toprule
Constraint on~$z=(x,y)$ & Effect on~$u$ & $m$ & Used in \\
\midrule
none
  & $u$ arbitrary
  & $k$
  & \cref{prop:single}, $G=\id$ \\
\addlinespace
$x_{k+j}=G(x_j)$
  & $u_{k+j}=G\circ u_j$
  & $2k$
  & \cref{prop:single} \\
\addlinespace
$x_{s+1}=y_s$
  & $u_{s+1}=T^su_1$
  & $t$
  & \cref{thm:trajectories}, $k=1$ \\
\addlinespace
$x_{ik+j}=G_i(x_j)$
  & $u_{ik+j}=G_i\circ u_j$
  & $k(N+1)$
  & \cref{thm:several,thm:pairwise} \\
\midrule
$x_1=y_1$
  & $u_1$ is a fixed point
  & $1$
  & not triangular \\
\bottomrule
\end{tabular}
\end{table}

\subsection{Preliminaries}\label{subsec:prelim}
Before proving the repair lemma, we collect some known definitions
and elementary remarks.
On~$[-1,1]^p$ we use the $\ell^1$ norm~$\abs{z}:=\sum_{i=1}^p \abs{z_i}$.
For vector-valued functions, the $L^1(\mu)$ norm is the integral of this $\ell^1$ norm. For
$\rho,\rho'\in\CM([-1,1]^p)$ the \emph{Wasserstein distance} is
\[
	d_W(\rho,\rho'):=\sup_h\abs{\int h\d\rho-\int h\d\rho'},
\]
where~$h$ ranges over the functions~$[-1,1]^p\to\R$ that are $1$-Lipschitz for the
$\ell^1$ norm. It metrizes weak convergence, and~$(\CM([-1,1]^p),d_W)$ is a compact
metric space~\cite[Chapter~6]{villani2009optimal}. For nonempty
sets~$P,Q\subseteq\CM([-1,1]^p)$, the \emph{Hausdorff distance} is
\[
	d_H(P,Q):=\max\Bigl\{\sup_{\rho\in P}\dist(\rho,Q),\ \sup_{\rho\in Q}\dist(\rho,P)\Bigr\},
	\qquad \dist(\rho,Q):=\inf_{\rho'\in Q}d_W(\rho,\rho').
\]
Note that~$d_H$ is a pseudometric, $d_H(P,Q)=0$ if and only if~$\overline P=\overline Q$,
and~$d_H$ is a metric on the nonempty closed sets, which form a compact space.

Two maps~$T,T' \in \CT$, on possibly different spaces,
are weakly equivalent if and only if~$d_H(P_k(T),P_k(T'))=0$ for every~$k$,
that is~$\overline{P_k(T)}=\overline{P_k(T')}$ for every~$k$.
A map~$T$ is a limit of~$(T_n)_{n\in \N}$ if and only if~$d_H(P_k(T_n),P_k(T))\to0$ for
every~$k$, that is~$\overline{P_k(T_n)}\to \overline{P_k(T)}$ for every~$k$.
Limits are unique exactly up to weak equivalence.

\begin{definition}[Modulus of continuity] \label{def:modulus}
A \emph{modulus of continuity}~$\omega\colon[0,\infty)\to[0,\infty)$ is a nondecreasing,
right-continuous function with~$\omega(0)=0$.
For a standard probability space~$(\Omega,\mu)$, let~$\CT_\omega(\mu)$ be the set of
maps~$T\colon\Linf{\mu}\to\Linf{\mu}$ such that
\[
	\norm{Tu-Tw}_{L^1(\mu)}\leq\omega(\norm{u-w}_{L^1(\mu)})\qquad\text{for all~$u,w\in\Linf{\mu}$},
\]
and let~$\CT_\omega$ be the collection of all such maps, over all standard
probability spaces.
\end{definition}

A sequence in~$\CT$
is \emph{$L^1$-equicontinuous} if it is contained in~$\CT_\omega$ for
a single~$\omega$.
Every~$T\in\CT$ lies in~$\CT_{\omega_T}$,
where~$\omega_T(s):=\sup\{\norm{Tu-Tw}_{L^1(\mu)}\colon
\norm{u-w}_{L^1(\mu)}\leq s\}$. Convexity of~$\Linf{\mu}$ makes~$\omega_T$
subadditive; since~$\omega_T(s)\to0$ as~$s\downarrow0$, it is right-continuous.

\begin{lemma} [Elementary] \label{lem:elementary}
The following facts are true.
\begin{enumerate}[(a)]
\item \textup{(Coupling.)} If $a,b\colon\Omega\to[-1,1]^p$ are measurable, then
\[
		d_W(\CD_\mu(a),\CD_\mu(b))\leq \norm{a-b}_{L^1(\mu)}.
\]
\item \textup{(Pointwise maps.)} Let~$g\colon[-1,1]^p\to\R$ be continuous. There is a
	modulus~$\omega_g$, depending only on~$g$, such that on every probability
	space~$(\Omega,\mu)$ and for all measurable~$a,b\colon\Omega\to[-1,1]^p$ we have
$
		\norm{g\circ a-g\circ b}_{L^1(\mu)}
			\leq\omega_g\p{\norm{a-b}_{L^1(\mu)}}.
$
\item \textup{(Push-forwards.)} Let~$\Psi\colon[-1,1]^p\to[-1,1]^q$ be continuous.
	Then~$\rho\mapsto\Psi_*\rho$ is uniformly continuous from~$\CM([-1,1]^p)$
	to~$\CM([-1,1]^q)$. In particular~$\overline{\Psi_*P}=\Psi_*\overline P$ for every
	$P\subseteq\CM([-1,1]^p)$, and~$d_H(\Psi_*P_n,\Psi_*P)\to0$
	whenever~$d_H(P_n,P)\to0$.
\end{enumerate}
\end{lemma}

\begin{proof}
(a) For every $1$-Lipschitz~$h$ we have
$\abs{\int h\circ a\d\mu-\int h\circ b\d\mu}\leq\int\abs{a-b}\d\mu$.
(b) Since~$[-1,1]^p$ is compact, the function~$g$ admits a continuity modulus~$\kappa_g$
and a maximum~$\norm{g}_\infty$. For any~$\delta>0$,
splitting~$\Omega$ according to whether~$\abs{a-b}\leq\delta$ and using
Markov's inequality, we get
$
	\norm{g\circ a-g\circ b}_{L^1(\mu)} \leq
		\kappa_g(\delta) + 2 \norm{g}_\infty \delta^{-1} \norm{a-b}_{L^1(\mu)}.
$
Hence~$\omega_g(s):=\inf_{\delta>0}\bigl(\kappa_g(\delta)+2\norm g_\infty s/\delta\bigr)$ works:
it is nondecreasing, vanishes at~$0$, and tends to~$0$ as~$s\to0$ (take~$\delta=\sqrt s$).
(c) The map~$\rho\mapsto\Psi_*\rho$ preserves weak convergence, so it is continuous
in the Wasserstein metric, hence uniformly continuous by compactness
of~$\CM([-1,1]^p)$.
\end{proof}

By~(b), for every continuous~$G\colon[-1,1]\to[-1,1]$ the pointwise
map~$u\mapsto G\circ u$ satisfies~$\norm{G\circ u-G\circ w}_{L^1(\mu)}\leq\omega_G(\norm{u-w}_{L^1(\mu)})$
on every probability space. Hence, for~$\Gamma$ as in \cref{prop:single}
and~$T\in\CT_\omega$, applying~(b) to~$F$ with~$a=(u,T(G\circ u))$ and~$b=(w,T(G\circ w))$ gives
\[
	\norm{\Gamma(T)u-\Gamma(T)w}_{L^1(\mu)}\leq\omega_F\bigl(s+\omega(\omega_G(s))\bigr),
	\qquad s:=\norm{u-w}_{L^1(\mu)},
\]
so~$\Gamma$ maps~$\CT$ into~$\CT$.

\begin{lemma}[Moduli and limits] \label{lem:moduli}
Let~$T_n\to T$, with~$T_n,T\in\CT$.
\begin{enumerate}[(a)]
\item If~$T_n\in\CT_\omega$ for all~$n$, for a modulus~$\omega$, then~$T\in\CT_\omega$.
\item The maps~$T_n$ and~$T$ have a common $L^1$ modulus of continuity.
\end{enumerate}
\end{lemma}

Part~(b) is a convenient criterion for equicontinuity: it suffices to check
that the limit is uniformly continuous.

\begin{proof}
(a) We use $2$-profiles. Fix an arbitrary~$\rho:=\CD_\mu(u,w,Tu,Tw)\in P_2(T)$.
There exists a sequence~$\rho_n:=\CD_{\mu_n}(u_n,w_n,T_n u_n,T_n w_n)\in P_2(T_n)$
with~$d_W(\rho_n,\rho)\to 0$, thus
\[
	\norm{u-w}_{L^1(\mu)}
	= \int \abs{x_1-x_2} \d\rho
	= \lim_{n\to\infty} \int \abs{x_1-x_2} \d\rho_n
	= \lim_{n\to\infty} \norm{u_n-w_n}_{L^1(\mu_n)}
\]
and similarly~$\norm{Tu-Tw}_{L^1(\mu)}=\lim_{n\to\infty}\norm{T_nu_n-T_nw_n}_{L^1(\mu_n)}$.
It remains to notice that for any convergent real sequences~$a_n\to a$ and~$b_n\to b$
with~$b_n\leq \omega(a_n)$, then~$b\leq \omega(a)$ by monotonicity
and right-continuity of~$\omega$.

(b) Let~$\omega$ be a continuous modulus for~$T$, and set
$\delta_n:=d_H(\overline{P_2(T_n)},\overline{P_2(T)})\to0$.
For~$\rho\in\CM([-1,1]^4)$, the functionals
\[
	a(\rho):=\int|x_1-x_2|\d\rho,
	\qquad b(\rho):=\int|y_1-y_2|\d\rho
\]
are~$1$-Lipschitz for~$d_W$. Every~$\rho\in\overline{P_2(T)}$ satisfies
$b(\rho)\leq\omega(a(\rho))$. Matching a profile generated by~$u,w$ under~$T_n$
with a law in~$\overline{P_2(T)}$ therefore gives
\[
	\norm{T_nu-T_nw}_{L^1(\mu_n)}
	\leq\omega\bigl(\norm{u-w}_{L^1(\mu_n)}+\delta_n\bigr)+\delta_n.
\]
Given an output tolerance, first take~$n$ large and then~$\norm{u-w}_{L^1(\mu_n)}$ small.
The remaining finitely many maps are individually uniformly continuous.
Hence~$\sup_n\omega_{T_n}(s)\to0$ as~$s\downarrow0$. The supremum of these
optimal moduli and~$\omega_T$ is subadditive, so it is a common modulus.
\end{proof}



\subsection{The repair lemma}
Let~$g$ be a triangular constraint, as in \cref{def:triangular}. The
set~$K_g\subseteq[-1,1]^{2m}$ is closed, and~$\varphi_g$ is a continuous
nonnegative function vanishing exactly on~$K_g$. For a
map~$T\colon\Linf{\mu}\to\Linf{\mu}$ and~$u\in\Linf{\mu}^m$, the
\emph{defect}
\begin{equation} \label{eq:defect}
	\eta_g(u):=\int\varphi_g\d\CD_\mu(u,Tu)
	=\sum_{j\in J}\norm{u_j-g_j\circ (u_{<j},Tu_{<j})}_{L^1(\mu)},
\end{equation}
where~$u_{<j}:=(u_1,\ldots,u_{j-1})$ and~$Tu_{<j}:=(Tu_1,\ldots,Tu_{j-1})$,
measures how far~$u$ is from satisfying the constraints.
In other words, we have~$\CD_\mu(u,Tu)\in P_m(T)\cap\CM(K_g)$
if and only if~$\eta_g(u)=0$, that is, if and only if
\begin{equation} \label{eq:constraint}
	u_j=g_j\circ(u_{<j},Tu_{<j})\quad\text{$\mu$-almost everywhere, for every }j\in J.
\end{equation}
The set~$P_m(T)\cap\CM(K_g)$ is nonempty: starting from any~$\CD_\mu(u,Tu)\in P_m(T)$,
replace~$u_j$ by~$g_j\circ(u_{<j},Tu_{<j})$, in increasing order of~$j\in J$.
The following key lemma shows that the cost of these consecutive repairs is controlled
by the initial defect~$\eta_g(u)$.

\begin{lemma}[Repair] \label{lem:repair}
Let~$\omega$ be a modulus and~$g$ be a triangular constraint on~$[-1,1]^{2m}$. Then there is a
modulus~$\Theta$, depending only on~$\omega$ and~$g$, with the following property.
For every probability space~$(\Omega,\mu)$, every~$T\in\CT_\omega(\mu)$ and
every~$u\in\Linf{\mu}^m$, there exists~$\tilde u\in\Linf{\mu}^m$ such
that~$\eta_g(\tilde u)=0$, that is,~$\CD_\mu(\tilde u,T\tilde u)\in P_m(T)\cap\CM(K_g)$, and
\[
	d_W\bigl(\CD_\mu(u,Tu),\CD_\mu(\tilde u,T\tilde u)\bigr)\leq\Theta\bigl(\eta_g(u)\bigr).
\]
\end{lemma}

\begin{proof}
For~$w,w'\in\Linf{\mu}^m$ let
\[
	\norm{w-w'}_T:=\sum_{i=1}^m\bigl(\norm{w_i-w'_i}_{L^1(\mu)}+\norm{Tw_i-Tw'_i}_{L^1(\mu)}\bigr)
\]
be the $L^1$-distance between~$(w,Tw)$ and~$(w',Tw')$. By \cref{lem:elementary}(a),
\begin{equation} \label{eq:dW_T}
	d_W\bigl(\CD_\mu(w,Tw),\CD_\mu(w',Tw')\bigr)\leq\norm{w-w'}_T.
\end{equation}
We fix two moduli once and for all: $\bar\omega(s):=s+\omega(s)$, and a
modulus~$\omega_{\varphi_g}$, depending only on~$g$, such that
$\abs{\eta_g(w)-\eta_g(w')}\leq\omega_{\varphi_g}(\norm{w-w'}_T)$ on every probability
space; it exists by \cref{lem:elementary}(b), applied to the continuous
function~$\varphi_g$.

Let~$w\in\Linf{\mu}^m$ and~$j\in J$, and let~$w'$ be
obtained from~$w$ by replacing its component~$w_j$
with~$g_j\circ(w_{<j},Tw_{<j})$. Since the components before~$j$ are
unchanged, $w'$ satisfies the constraint at~$j$, and every constraint at~$i<j$
that~$w$ already satisfies. The key observation is that the size of the repair is
bounded by the defect: by~\eqref{eq:defect},
\[
	\norm{w_j-w'_j}_{L^1(\mu)}=\norm{w_j-g_j\circ(w_{<j},Tw_{<j})}_{L^1(\mu)}\leq\eta_g(w).
\]
Since only the components~$w_j$ and~$Tw_j$ change,
\[
	\norm{w-w'}_T=\norm{w_j-w'_j}_{L^1(\mu)}+\norm{Tw_j-Tw'_j}_{L^1(\mu)},
\]
hence, by definition of~$\bar\omega$ and~$\omega_{\varphi_g}$,
\[
	\norm{w-w'}_T\leq\bar\omega\bigl(\eta_g(w)\bigr),
	\qquad
	\eta_g(w')\leq\eta_g(w)+\omega_{\varphi_g}\bigl(\bar\omega(\eta_g(w))\bigr)
	=:\theta\bigl(\eta_g(w)\bigr).
\]

Let~$J=\{j_1<\cdots<j_r\}$, and let~$u=u^0,u^1,\ldots,u^r=:\tilde u$,
where~$u^\ell$ is obtained from~$u^{\ell-1}$ by repairing the component~$j_\ell$.
Then~$u^\ell$ satisfies the first~$\ell$ constraints exactly; in
particular~$\eta_g(\tilde u)=0$, which is the first claim.
Moreover~$\eta_g(u^\ell)\leq\theta^\ell(\eta_g(u))$, where~$\theta^\ell$ is the
$\ell$-fold composition, so by the triangle inequality
\[
	\norm{u-\tilde u}_T\leq\sum_{\ell=0}^{r-1}\bar\omega\bigl(\theta^\ell(\eta_g(u))\bigr)
	=:\Theta\bigl(\eta_g(u)\bigr).
\]
The function~$\Theta$ is a finite sum of compositions of the fixed moduli~$\omega$
and~$\omega_{\varphi_g}$, hence a modulus depending only on~$\omega$ and~$g$.
By~\eqref{eq:dW_T}, this proves the second claim.
\end{proof}

The modulus~$\omega$ enters the proof at a single point, to pass from a small change
of an input to a small change of the corresponding output; since~$\Theta$ depends
on~$T$ only through~$\omega$, the repair is uniform over an equicontinuous family.
We write
\[
	P^g_m(T):=P_m(T)\cap\CM(K_g).
\]

\begin{proposition}[Continuity of constrained profiles] \label{prop:constrained}
Let~$\omega$ be a modulus and~$g$ a triangular constraint on~$[-1,1]^{2m}$. Then there is a
modulus~$\omega'$, depending only on~$\omega$ and~$g$, such that for
all~$T\in\CT_\omega(\mu)$ and~$T'\in\CT_\omega(\mu')$,
\[
	d_H\bigl(P^g_m(T),P^g_m(T')\bigr)\leq\omega'\bigl(d_H(P_m(T),P_m(T'))\bigr).
\]
\end{proposition}

\begin{proof}
Let~$\kappa$ be a modulus of uniform continuity of~$\rho\mapsto\int\varphi_g\d\rho$
on~$\CM([-1,1]^{2m})$. Let~$\Theta$ be given by \cref{lem:repair}.
Let~$d:=d_H(P_m(T),P_m(T'))$, let~$s>d$, and let~$\rho'\in P^g_m(T')$. Then there
exists~$\rho\in P_m(T)$ with~$d_W(\rho,\rho')<s$. By \cref{lem:repair} there
exists~$\tilde\rho\in P^g_m(T)$
with $d_W(\tilde\rho,\rho)\leq \Theta\p{\int \varphi_g \d\rho}$.
Since~$\int \varphi_g \d\rho\leq \kappa(d_W(\rho,\rho'))\leq\kappa(s)$
as~$\int \varphi_g \d\rho'=0$, we conclude that
$
	d_W(\tilde\rho,\rho') \leq s + \Theta(\kappa(s)).
$
Letting~$\omega'(d):=\inf_{s>d}\bigl(s+\Theta(\kappa(s))\bigr)$, we
get~$\dist(\rho',P^g_m(T))\leq\omega'(d)$ for every~$\rho'\in P^g_m(T')$.
Exchanging~$T$ and~$T'$ completes the proof.
\end{proof}

\begin{corollary}[Lifting] \label{cor:lifting}
Let~$(T_n)_{n\in \N}$ be an $L^1$-equicontinuous sequence action convergent to~$T$. Let~$g$ be
a triangular constraint and~$\Psi\colon[-1,1]^{2m}\to[-1,1]^q$ continuous. Then
\[
	d_H\bigl(\Psi_*P^g_m(T_n),\Psi_*P^g_m(T)\bigr)\to0 .
\]
\end{corollary}

\begin{proof}
By \cref{lem:moduli}(a), the maps~$T_n$ and~$T$ have a common modulus.
Combine \cref{prop:constrained} and \cref{lem:elementary}(c).
\end{proof}


\subsection{Convergence of trajectories}
As a first application of the lifting trick, we prove \cref{thm:trajectories}, on the
joint distributions of finite trajectories.

\begin{proof}[Proof of \cref{thm:trajectories}]
Fix~$t,k\geq1$, and index~$m=kt$ coordinates by pairs~$(s,j)$, with~$1\leq s\leq t$
and~$1\leq j\leq k$, ordered lexicographically with~$s$ first. Consider the triangular
constraint given by~$x_{(s+1,j)}=y_{(s,j)}$ for~$1\leq s<t$ and all~$j$:
the next input is the previous output. Let~$\Psi$ be the coordinate projection
\[
	\Psi(x,y):=\bigl((x_{(1,j)})_{j},(y_{(1,j)})_{j},(y_{(2,j)})_{j},\ldots,(y_{(t,j)})_{j}\bigr).
\]
By induction on~$s$, a configuration~$w\in\Linf{\mu}^{m}$
satisfies~\eqref{eq:constraint} if and only if~$w_{(s+1,j)}=T^sw_{(1,j)}$ for all~$s<t$
and~$j$, and then~$\Psi_*\CD_\mu(w,Tw)=\CD_\mu(u,Tu,\ldots,T^tu)$
with~$u:=(w_{(1,j)})_{j\leq k}$; every~$u\in\Linf{\mu}^k$ arises in this way.
Hence~$P^{(t)}_k(T)=\Psi_*P^g_m(T)$ for every map~$T$, and \cref{cor:lifting} applies.
\end{proof}


\subsection{Simple couplings}\label{subsec:composition_proof}
We first treat the case in which the coupling function depends only on the
neighbour's state. Here linearity is not needed.
The following theorem already covers
coupled systems~\eqref{eq:intro_nonlin}
in the case~$G(x,y)=h(y)$.

\begin{theorem}  \label{prop:single}
Let~$F\colon[-1,1]^2\to[-1,1]$ and~$G\colon[-1,1]\to[-1,1]$ be continuous.
For~$T\in\CT$, define~$\Gamma(T)\colon\Linf{\mu}\to\Linf{\mu}$ by
\[
	(\Gamma(T)u)(x):=F\bigl(u(x),\,T(G\circ u)(x)\bigr).
\]
Then~$\Gamma(T)\in\CT$. Moreover,
if~$(T_n)$ is $L^1$-equicontinuous and action convergent to~$T$,
then~$(\Gamma(T_n))$ is action convergent to~$\Gamma(T)$.
\end{theorem}

\begin{proof}
We obtain convergence of the $k$-profile of~$\Gamma(T_n)$ from convergence
of the $2k$-profile of~$T_n$. Fix any~$T\in\CT$.
Let~$g$ be the constraint on~$m=2k$ coordinates
given by the equations~$x_{k+j}=G(x_j)$ for~$1\leq j\leq k$,
so that a point~$\CD_\mu(w,Tw)\in P_{2k}(T)$ belongs to~$\CM(K_g)$ if and only
if~$w=(u,G\circ u)$ for some~$u\in\Linf{\mu}^k$, with~$G$ acting coordinatewise. Let
$
	\Psi(x,y):=\bigl(x_1,\ldots,x_k,\,F(x_1,y_{k+1}),\ldots,F(x_k,y_{2k})\bigr).
$
Then~$\Psi_*\CD_\mu(w,Tw)=\CD_\mu(u,\Gamma(T)u)$, and every~$u$ arises in this way.
Hence~$P_k(\Gamma(T))=\Psi_*P^g_{2k}(T)$. This holds for any map~$T$, and
\cref{cor:lifting} applies.
\end{proof}

\begin{remark}
If~$G=\id$, let~$\Psi(x,y):=\bigl(x,F(x_1,y_1),\ldots,F(x_k,y_k)\bigr)$.
Then~$P_k(\Gamma(T))=\Psi_*P_k(T)$ for every map~$T$, and the
claim follows from \cref{lem:elementary}(c), without equicontinuity.
For general~$G$, equicontinuity cannot be dropped.
\end{remark}

The proof applies verbatim to several inner functions:

\begin{theorem} \label{thm:several}
Let
$G_1,\ldots,G_N\colon[-1,1]\to[-1,1]$ and~$F\colon[-1,1]^{1+N}\to[-1,1]$
be continuous. For any~$T\colon\Linf{\mu}\to\Linf{\mu}$ in~$\CT$,
define~$\Gamma(T)\colon\Linf{\mu}\to\Linf{\mu}$ by
\[
	(\Gamma(T)u)(x):=F\bigl(u(x),\,T(G_1\circ u)(x),\ldots,T(G_N\circ u)(x)\bigr).
\]
Then~$\Gamma$ maps~$\CT$ into~$\CT$ and preserves action convergence of
$L^1$-equicontinuous sequences.
\end{theorem}

\begin{proof}
The proof is similar, using the constraint~$x_{ik+j}=G_i(x_j)$
for~$1\leq i\leq N$, $1\leq j\leq k$ on~$m=k(N+1)$ coordinates,
and~$\Psi\colon[-1,1]^{2m}\to[-1,1]^{2k}$ defined, for~$1\leq j\leq k$,
by~$\Psi(x,y)_j:=x_j$
and~$\Psi(x,y)_{k+j}:=F\bigl(x_j,\,y_{k+j},\,y_{2k+j},\ldots,y_{Nk+j}\bigr)$.
\end{proof}


\subsection{Pairwise couplings}\label{subsec:pairwise}
We now prove \cref{thm:pairwise}.
In contrast to \cref{prop:single}, the underlying transformation is now required
to be linear.
Recall the class~$\CL$ of linear maps in~$\CT$, from \cref{subsec:definitions}.
The operators of \cref{ex:dense} and \cref{ex:randomwalk}, the averaging operator~$J$,
and all Markov operators preserving~$\mu$ belong to~$\CL$, with~$\omega(s)=s$.

We now define, for~$W\in\CL$ and continuous~$G$, an
operator~$W[G]\colon\Linf{\mu}\to L^\infty(\mu)$ playing the role of the inner
sum in~\eqref{eq:network}.
For a \emph{tensor}~$H=\sum_la_l\otimes b_l$,
with~$(a\otimes b)(x,y):=a(x)b(y)$
a finite sum with~$a_l,b_l$
continuous and real-valued, define
\[
	W[H]u:=\sum_l(a_l\circ u)\cdot W(b_l\circ u)\in L^\infty(\mu).
\]

To extend~$H\mapsto W[H]$ to all continuous~$G$, we need a bound
uniform in~$u$.

\begin{lemma}[Tensor bound] \label{lem:tensor}
Let~$W$ be a linear operator on~$L^\infty(\mu)$ with~$\norm{Wu}_\infty\leq\norm u_\infty$.
Then~$\abs{W[H]u}\leq\norm H_\infty$ $\mu$-almost everywhere, for all
tensors~$H$ and all~$u\in\Linf{\mu}$. In particular, $W[H]$ depends only
on the function~$H$, not on its representation.
\end{lemma}

\begin{proof}
If~$u=\sum_pc_p\1_{E_p}$ is simple, then on~$E_p$, by linearity,
$W[H]u=W\bigl(H(c_p,u(\cdot))\bigr)$, which is bounded
by~$\norm{H(c_p,u(\cdot))}_\infty\leq\norm H_\infty$.
For general~$u$, approximate it uniformly by simple functions~$u_r$: then
$a_l\circ u_r\to a_l\circ u$ and~$b_l\circ u_r\to b_l\circ u$ uniformly,
hence~$W(b_l\circ u_r)\to W(b_l\circ u)$ in~$L^\infty$, so~$W[H]u_r\to W[H]u$
in~$L^\infty$ and the bound passes to the limit.
For the last claim, apply the bound to the difference of two representations.
\end{proof}

By the Stone--Weierstrass theorem, tensors are dense in~$C([-1,1]^2)$. Moreover, every
tensor can be written with all~$a_l,b_l$ taking values in~$[-1,1]$: for an
integer~$K\geq\norm{a_l}_\infty$, replace~$a_l\otimes b_l$ by~$K$ copies
of~$(a_l/K)\otimes b_l$, and similarly for~$b_l$. Since~$H\mapsto W[H]u$ is linear
on tensors, \cref{lem:tensor} shows that, for each~$u$, it extends uniquely to a
linear map~$C([-1,1]^2)\to L^\infty(\mu)$, and
\begin{equation} \label{eq:tensor_error}
	\norm{W[G]u-W[H]u}_\infty\leq\norm{G-H}_\infty
	\quad\text{for continuous~$G,H$ and~$u\in\Linf{\mu}$}.
\end{equation}
Thus~$W[H]\to W[G]$ uniformly on~$\Linf{\mu}$ as~$H\to G$ uniformly.
In particular, $W[G]$ maps~$\Linf{\mu}$ into itself if~$\abs G\leq1$.
If~$H(a,b)=h(b)$, then~$W[H]u=W(h\circ u)$; in particular~$W[H]=W$ for~$H(a,b)=b$.

\begin{remark}[Why linearity] \label{rem:linearity}
One might try to define pairwise coupling directly,
as a map~$x\mapsto W\bigl(G(u(x),u(\cdot))\bigr)(x)$. This is not well defined,
not even~$W=\id$: the input~$G(u(x),u(\cdot))$ depends on~$x$, and the output is
evaluated at the same point~$x$, so the formula uses values of elements
of~$L^\infty(\mu)$ along the diagonal of~$\Omega\times\Omega$, which may be a null set.
For tensors, linearity moves~$x$ outside~$W$, and definition follows by density.
\end{remark}

While not needed to prove convergence, it is convenient to represent~$W\in \CL$
as an integral operator, as we did in \cref{thm:pairwise} to simplify exposition.
Such representation is given by the following lemma.

\begin{lemma}[Signed kernels]\label{lem:signed_kernel}
Every~$W\in\CL$ admits a measurable family~$(\nu_x)_{x\in\Omega}$ of
finite signed measures with~$|\nu_x|(\Omega)\leq1$ such that~$Wv(x)=\int_\Omega v\d\nu_x$
for~$\mu$-almost every~$x$, for every~$v\in L^\infty(\mu)$. Moreover,
$\mu(B)=0$ implies~$|\nu_x|(B)=0$ for~$\mu$-almost every~$x$.
\end{lemma}

\begin{proof}
This is a special case of the representation of operators by random measures,
due to Sourour; see~\cite[Theorem~2.1 and Remarks~2.2, 2.3]{weis1984representation}.
\end{proof}

The class~$\CL$ is closed under limits of equicontinuous sequences:

\begin{lemma}\label{lem:linear_limit}
Let~$W_n\in\CL$ have a common modulus~$\omega$. If~$W_n\to W$,
then~$W\in\CL\cap\CT_\omega$.
\end{lemma}

\begin{proof}
By \cref{lem:moduli}(a), $W\in\CT_\omega$, and by \cref{prop:constrained},
$d_H(P^g_m(W_n),P^g_m(W))\to0$ for every triangular constraint~$g$. Take the
constraint~$x_3=(x_1+x_2)/2$ on~$m=3$ coordinates. By linearity of~$W_n$, every law
in~$P^g_3(W_n)$ is supported on the closed set~$\{y_3=(y_1+y_2)/2\}$, hence so is
every law in~$P^g_3(W)$; that is,~$W\bigl(\frac{u+w}2\bigr)=\frac{Wu+Ww}2$ for
all~$u,w\in\Linf{\mu}$. In the same way, the constraint~$x_2=cx_1$
gives~$W(cu)=cWu$ for~$c\in[-1,1]$. Hence~$Lv:=\norm v_\infty W(v/\norm v_\infty)$,
with~$L0:=0$, defines a linear operator on~$L^\infty(\mu)$ extending~$W$,
with~$\norm{Lv}_\infty\leq\norm v_\infty$ since~$W$ maps~$\Linf{\mu}$ into itself.
Hence~$W\in\CL$.
\end{proof}

\begin{example}
For~$W\in\CL$ represented as in \cref{lem:signed_kernel}, we have
\begin{equation} \label{eq:kernel}
	W[G]u(x)=\int G\bigl(u(x),u(y)\bigr)\d\nu_x(y).
\end{equation}
The identity holds first for tensors and then for continuous~$G$ by uniform
approximation and the bound~$|\nu_x|(\Omega)\leq1$.
For the constant kernel~$J$ it is~$\int G(u(x),u(y))\d\lambda(y)$, the limit
used in \cref{cor:STO}. In particular, the map~$\Gamma(W)$ of \cref{thm:pairwise}
is~$u\mapsto F(u,W[G]u)$.
\end{example}

We are now ready to prove \cref{thm:pairwise}.

\begin{proof}[Proof of \cref{thm:pairwise}]
By \cref{thm:existence}, the sequence~$(W_n)$ has a limit~$W\in\CT$, and by
\cref{lem:linear_limit} it is again in~$\CL$, with the
common modulus~$\omega$ of the sequence; by \eqref{eq:kernel}, $\Gamma(W)u=F(u,W[G]u)$.
The idea is to replace~$G$ by a tensor, which brings us to \cref{thm:several}.
For convenience, extend~$F$ to~$[-1,1]\times\R$ by~$F(a,c):=F\bigl(a,\max\{-1,\min\{1,c\}\}\bigr)$;
the extension is still continuous, and a modulus~$\kappa_F$ with
$\abs{F(a,c)-F(a,c')}\leq\kappa_F(\abs{c-c'})$ for~$c,c'\in[-1,1]$ works for
all~$c,c'\in\R$, since the truncation is $1$-Lipschitz.
Let~$\ve>0$, and take a tensor~$H=\sum_{l=1}^Na_l\otimes b_l$
with~$\norm{G-H}_\infty\leq\ve$ and all~$a_l,b_l$ taking values in~$[-1,1]$.
Let~$\Gamma_\ve$ be the map of \cref{thm:several} with outer function
$
	\Phi(a,c_1,\ldots,c_N):=F\Bigl(a,\ \sum_{l=1}^Na_l(a)\,c_l\Bigr)
$
and inner functions~$b_1,\ldots,b_N$, so that~$\Gamma_\ve(W)u=F\bigl(u,W[H]u\bigr)$.
By~\eqref{eq:tensor_error} and the choice of~$\kappa_F$,
\begin{equation} \label{eq:pairwise_approx}
	\norm{\Gamma_\ve(W)u-\Gamma(W)u}_\infty\leq\kappa_F(\ve)
	\qquad\text{for all~$W\in\CL$ and~$u\in\Linf{\mu}$}.
\end{equation}

\emph{$\Gamma(W)\in\CT$.} By \cref{thm:several}, $\Gamma_\ve(W)$ has a
modulus~$\omega_\ve$ depending only on~$F,G,\ve$ and a modulus for~$W$, so~\eqref{eq:pairwise_approx}
gives
\[
\norm{\Gamma(W)u-\Gamma(W)w}_{L^1(\mu)}\leq2\kappa_F(\ve)+\omega_\ve(\norm{u-w}_{L^1(\mu)}).
\]
The infimum over~$\ve>0$ of the right-hand side is a modulus, since an infimum of
nondecreasing right-continuous functions is again right-continuous.

Fix~$k\geq1$. By~\eqref{eq:pairwise_approx} and \cref{lem:elementary}(a),
$d_H(P_k(\Gamma_\ve(W)),P_k(\Gamma(W)))\leq k\kappa_F(\ve)$ for every~$W\in\CL$.
The sequence~$(W_n)$ is $L^1$-equicontinuous, so
\cref{thm:several} gives
\[
	d_H(P_k(\Gamma_\ve(W_n)),P_k(\Gamma_\ve(W)))\to0.
\]
By the triangle inequality
\[
	\limsup_nd_H(P_k(\Gamma(W_n)),P_k(\Gamma(W)))\leq2k\kappa_F(\ve),
\]
and~$\ve$ is arbitrary.
\end{proof}


\subsection{A converse to \texorpdfstring{\cref{thm:pairwise}}{Theorem B}}\label{subsec:converse}
Under a mild nondegeneracy assumption, the limit of a single coupled system
determines the limit of the underlying linear operators.
The proof is again a ``lift and repair'' argument.

\begin{proposition}[Converse of \cref{thm:pairwise}] \label{prop:converse}
Let~$F,G\colon[-1,1]^2\to[-1,1]$ be continuous, and suppose that there
are~$a_1\neq a_2$ such that, for~$i=1,2$, we have~$G(a_i,a_1)\neq G(a_i,a_2)$ and
$c\mapsto F(a_i,c)$ is injective. Let~$W_n,W\in\CL$ with~$W_n\1=W\1=\1$.
Then~$\Gamma(W_n)\to\Gamma(W)$ if and only if~$W_n\to W$. In particular, if~$W'\in\CL$
satisfies~$W'\1=\1$ and~$\Gamma(W')$ is weakly equivalent to~$\Gamma(W)$, then~$W'$ is
weakly equivalent to~$W$.
\end{proposition}

\begin{proof}
If~$W_n\to W$, then~$\Gamma(W_n)\to\Gamma(W)$ by \cref{thm:pairwise}. Conversely,
let~$\Gamma(W_n)\to\Gamma(W)$. By \cref{thm:pairwise}, these maps lie in~$\CT$, hence
in~$\CT_\omega$ for a single modulus~$\omega$, by \cref{lem:moduli}(b).
We write the proof for~$k=1$; the case~$k>1$ is similar.

The idea is to read~$P_1(V)$, for~$(\Omega,\mu,V)\in\CL$ with~$V\1=\1$, from the laws
in~$P_N(\Gamma(V))$ whose inputs take values in~$\{a_1,a_2\}$, through a continuous
map~$\Phi$ that does not depend on~$V$. Applied to~$V=W_n$ and~$V=W$, this reduces the
claim to the convergence of these restricted profiles, which follows by a repair.

Let~$\Delta_i:=G(a_i,a_2)-G(a_i,a_1)\neq0$.
Since~$F(a_i,\cdot)$ is a homeomorphism of~$[-1,1]$ onto a compact interval, there is a
continuous~$\psi\colon\{a_1,a_2\}\times[-1,1]\to\R$ with
$\psi(a_i,F(a_i,c))=(c-G(a_i,a_1))/\Delta_i$ for all~$c\in[-1,1]$.
For measurable~$E\subseteq\Omega$
let~$u_E:=a_1\1_{\Omega\setminus E} + a_2 \1_{E}$.
From~\eqref{eq:kernel}, $V\1=\1$ (which gives~$\nu_x(\Omega)=1$),
and the definition of~$\psi$, we have
\[
	V\1_E (x) = \nu_x(E)
	= \psi\bigl(u_E(x),\Gamma(V)u_E(x)\bigr)\qquad\text{$\mu$-almost every~$x\in \Omega$.}
\]

Fix~$N\geq 1$. Let~$\clip(s):=\max\{-1,\min\{1,s\}\}$.
On~$K:=\{a_1,a_2\}^N\times[-1,1]^N$, define the
continuous map~$\Phi\colon K\to[-1,1]^2$ by
\[
	\Phi(x,y):=\Bigl(-1+\frac2N\sum_{j=1}^N\1_{\{x_j=a_2\}},\
	\clip\Bigl(-1+\frac2N\sum_{j=1}^N\psi(x_j,y_j)\Bigr)\Bigr).
\]
Given measurable~$E_1,\ldots,E_N\subseteq\Omega$, let~$u:=(u_{E_1},\ldots,u_{E_N})$
and~$v:=-1+\frac2N\sum_j\1_{E_j}$, which lies in~$\Linf{\mu}$. By linearity
and~$V\1=\1$, we have~$Vv=-1+\frac2N\sum_jV\1_{E_j}$, so the previous step
gives~$\Phi_*\CD_\mu(u,\Gamma(V)u)=\CD_\mu(v,Vv)$. Every~$u$
with~$\CD_\mu(u,\Gamma(V)u)\in\CM(K)$ is of this form, hence
$R(V):=\Phi_*\bigl(P_N(\Gamma(V))\cap\CM(K)\bigr)\subseteq P_1(V)$.
Conversely, for~$w\in\Linf{\mu}$, the sets~$E_j:=\{w>-1+2j/N\}$
give~$\norm{w-v}_\infty\leq2/N$, hence~$\norm{Vw-Vv}_\infty\leq2/N$.
By \cref{lem:elementary}(a), we have~$d_H(P_1(V),R(V))\leq4/N$.

By the previous step, applied to~$V=W_n$ and~$V=W$,
\[
	d_H\bigl(P_1(W_n),P_1(W)\bigr)\leq\frac8N+d_H\bigl(R(W_n),R(W)\bigr).
\]
Write~$P^K(T):=P_N(T)\cap\CM(K)$, so that~$R(V)=\Phi_*P^K(\Gamma(V))$, and
let~$\kappa$ be a nondecreasing modulus of uniform continuity of~$\rho\mapsto\Phi_*\rho$
on~$\CM(K)$ (\cref{lem:elementary}(c)). It remains to show that the restricted
profiles~$P^K(\Gamma(W_n))$ and~$P^K(\Gamma(W))$ are close. More precisely, we claim that
for every~$s>d_H\bigl(P_N(\Gamma(W_n)),P_N(\Gamma(W))\bigr)$,
\[
	d_H\bigl(P^K(\Gamma(W_n)),P^K(\Gamma(W))\bigr)\leq2s+N\omega(s).
\]
Let~$\rho'\in P^K(\Gamma(W))$. There is~$\rho=\CD_{\mu_n}(u,\Gamma(W_n)u)$
with~$d_W(\rho,\rho')<s$. The function~$\varphi(x,y):=\sum_j\dist(x_j,\{a_1,a_2\})$ is
$1$-Lipschitz and vanishes on~$K$, so~$\eta:=\int\varphi\d\rho$ satisfies~$\eta<s$.
We now repair~$u$:
let~$\tilde u$ be obtained by replacing each value of~$u$ with a nearest point
of~$\{a_1,a_2\}$. Then~$\tilde\rho:=\CD_{\mu_n}(\tilde u,\Gamma(W_n)\tilde u)\in P^K(\Gamma(W_n))$
and~$\sum_j\norm{u_j-\tilde u_j}_{L^1(\mu_n)}=\eta$, so by \cref{lem:elementary}(a)
$
	d_W(\tilde\rho,\rho')<s+\eta+N\omega(\eta)\leq2s+N\omega(s).
$
The same holds with~$W_n$ and~$W$ interchanged, which proves the claim.
Combining the two displays,
\[
	d_H\bigl(P_1(W_n),P_1(W)\bigr)\leq\frac8N+\kappa\bigl(2s+N\omega(s)\bigr).
\]
Since~$\Gamma(W_n)\to\Gamma(W)$, by the definition of~$s$ we can let~$s\to0$ as~$n\to\infty$,
thus obtaining~$\limsup_nd_H(P_1(W_n),P_1(W))\leq8/N$ for every~$N$,
then we let~$N\to \infty$. The last claim is the case~$W_n=W'$.
\end{proof}


\section{Proof of \texorpdfstring{\cref{thm:existence}}{Theorem A}}\label{sec:existence}

Throughout this section, $(T_n)_{n\in\N}$ is an action convergent sequence
with~$T_n\in\CT_\omega(\mu_n)$ for a fixed modulus~$\omega$, where~$T_n$ acts
on~$\Linf{\mu_n}$ over~$(\Omega_n,\mu_n)$. For every~$k\geq1$, let~$Q_k$ be the
limit of the closed sets~$\overline{P_k(T_n)}$. It is a nonempty closed subset
of~$\CM([-1,1]^{2k})$, and~$d_H(P_k(T_n),Q_k)\to0$. We construct a standard probability
space~$(\Omega,\mu)$ and a map~$T\in\CT_\omega(\mu)$ with~$\overline{P_k(T)}=Q_k$ for
every~$k$. That every limit lies in~$\CT_\omega$ is \cref{lem:moduli}(a), and the
converse implication of \cref{thm:existence} is \cref{lem:moduli}(b).

We adapt the construction of Backhausz and Szegedy~\cite{backhausz2022action}.
On each~$\Omega_n$ we fix a countable family of functions whose input--output laws
approximate a dense subset of every~$Q_k$ as~$n\to\infty$, and which is closed
under~$T_n$ and under clipped polynomials with rational coefficients. Their joint laws live on a compact space, so they have a limit
point~$\mu$. On the limit space the family becomes the family of coordinate functions:
closure under~$T_n$ defines~$T$ on the coordinates, and closure under polynomials makes
the coordinates dense, so that~$T$ extends to the whole ball.

\emph{Words.}
Let~$\clip(s):=\max\{-1,\min\{1,s\}\}$.
For every~$k\geq 1$, fix a countable dense subset~$Q'_k\subseteq Q_k$.
Let~$\CW$ be the smallest set of formal expressions, called \emph{words}, such that:
\begin{enumerate}[(W1)]
\item $(t,j)\in\CW$ for all~$k\geq1$, $t\in Q'_k$ and~$1\leq j\leq k$;
\item $T\alpha\in\CW$ for every~$\alpha\in\CW$;
\item $p(\alpha_1,\ldots,\alpha_r)\in\CW$ for all~$r\geq1$,
	$\alpha_1,\ldots,\alpha_r\in\CW$ and~$p\in\Q[z_1,\ldots,z_r]$.
\end{enumerate}
Here~$T\alpha$ and~$p(\alpha_1,\ldots,\alpha_r)$ are formal symbols, distinct for
distinct data; the letter~$T$ anticipates the map constructed below. The set~$\CW$
is countable by construction.

The main difference from~\cite{backhausz2022action} is the choice of~$\CW$.
Linear maps commute with linear combinations, so in~\cite{backhausz2022action}
it suffices to close the words under products. Nonlinear maps do not, so we
close the words under linear combinations as well, that is, under polynomials.

For any~$n\geq 1$ and~$\alpha\in \CW$, we
define~$v_{n,\alpha}\in\Linf{\mu_n}$
by recursion on the construction of~$\alpha$. For~$t\in Q'_k$, choose~$u\in\Linf{\mu_n}^k$
with~$d_W(\CD_{\mu_n}(u,T_nu),t)\leq d_H(P_k(T_n),Q_k)+1/n$ and
set~$v_{n,(t,j)}:=u_j$ for~$1\leq j\leq k$, so that
\begin{equation} \label{eq:generators}
	\CD_{\mu_n}\bigl((v_{n,(t,j)})_{j\leq k},(T_nv_{n,(t,j)})_{j\leq k}\bigr)\to t
	\qquad\text{as~$n\to\infty$}.
\end{equation}
Then set~$v_{n,T\alpha}:=T_nv_{n,\alpha}$
and~$v_{n,p(\alpha_1,\ldots,\alpha_r)}
	:=\clip\circ p\circ (v_{n,\alpha_1},\ldots,v_{n,\alpha_r})$.
To shorten notation, let~$\tilde p:=\clip\circ p$.
The clipping keeps the values in~$[-1,1]$.

\emph{The limit space.} Let~$\Omega:=[-1,1]^{\CW}$ with the product topology,
and let~$\pi_\alpha\colon\Omega\to[-1,1]$ denote the coordinate of index~$\alpha$.
The Borel $\sigma$-algebra of~$\Omega$ is generated by the coordinates,
thus~$\Omega_n\to\Omega$, $x\mapsto(v_{n,\alpha}(x))_{\alpha\in\CW}$ is a measurable map.
Let~$\kappa_n$ be its law.
Its image under~$(\pi_{\alpha_1},\ldots,\pi_{\alpha_r})$ is the
law~$\CD_{\mu_n}(v_{n,\alpha_1},\ldots,v_{n,\alpha_r})$.
By compactness of~$\CM(\Omega)$,
there exists a subsequence~$(n_i)_{i\in\N}$ and a Borel probability measure~$\mu$
on~$\Omega$ with~$\int f\d\kappa_{n_i}\to\int f\d\mu$ for every continuous~$f\colon\Omega\to\R$.
If~$\Psi\colon\Omega\to[-1,1]^p$
is continuous, then~$\Psi_*\kappa_{n_i}\to\Psi_*\mu$ weakly, hence in~$d_W$.

The following two lemmas are analogous to~\cite[Lemma 4.2]{backhausz2022action}.

\begin{lemma}[Limit law] \label{lem:limit_law}
Let~$\alpha,\beta,\alpha_1,\ldots,\alpha_r\in\CW$ and~$p\in\Q[z_1,\ldots,z_r]$. Then:
\begin{enumerate}[(a)]
\item $\pi_{p(\alpha_1,\ldots,\alpha_r)}=\tilde p\circ (\pi_{\alpha_1},\ldots,\pi_{\alpha_r})$
	$\mu$-almost everywhere;
\item $\norm{\pi_{T\alpha}-\pi_{T\beta}}_{L^1(\mu)}\leq\omega\bigl(\norm{\pi_\alpha-\pi_\beta}_{L^1(\mu)}\bigr)$;
\item $\CD_\mu\bigl((\pi_{\alpha_j})_{j\leq k},(\pi_{T\alpha_j})_{j\leq k}\bigr)\in Q_k$
	for all~$\alpha_1,\ldots,\alpha_k\in\CW$;
\item $\CD_\mu\bigl((\pi_{(t,j)})_{j\leq k},(\pi_{T(t,j)})_{j\leq k}\bigr)=t$
	for all~$k\geq1$ and~$t\in Q'_k$.
\end{enumerate}
\end{lemma}

\begin{proof}
Each claim involves finitely many coordinates, through continuous functions on~$\Omega$.
(a) The continuous
function~$f:=\abs{\pi_{p(\alpha_1,\ldots,\alpha_r)}-\tilde p \circ (\pi_{\alpha_1},\ldots,\pi_{\alpha_r})}$
satisfies~$\int f\d\kappa_n=0$ for all~$n$, by definition of~$v_{n,p(\alpha_1,\ldots,\alpha_r)}$.
Hence~$\int f\d\mu=0$.
(b) Let~$a(\kappa):=\int\abs{\pi_\alpha-\pi_\beta}\d\kappa$
and~$b(\kappa):=\int\abs{\pi_{T\alpha}-\pi_{T\beta}}\d\kappa$. By definition
of~$v_{n,T\alpha}$ and since~$T_n\in\CT_\omega(\mu_n)$,
\[
	b(\kappa_n)=\norm{T_nv_{n,\alpha}-T_nv_{n,\beta}}_{L^1(\mu_n)}
	\leq\omega\bigl(\norm{v_{n,\alpha}-v_{n,\beta}}_{L^1(\mu_n)}\bigr)=\omega(a(\kappa_n)).
\]
Since~$a(\kappa_{n_i})\to a(\mu)$ and~$b(\kappa_{n_i})\to b(\mu)$, we get~$b(\mu)\leq\omega(a(\mu))$.
(c) Let~$\Psi(z):=\bigl((z_{\alpha_j})_{j\leq k},(z_{T\alpha_j})_{j\leq k}\bigr)$, a
continuous map~$\Omega\to[-1,1]^{2k}$. Then
\[
	\Psi_*\kappa_n=\CD_{\mu_n}\bigl((v_{n,\alpha_j})_{j\leq k},(T_nv_{n,\alpha_j})_{j\leq k}\bigr)
	\in P_k(T_n),
\]
and~$\Psi_*\kappa_{n_i}\to\Psi_*\mu$ in~$d_W$. Since~$d_H(P_k(T_{n_i}),Q_k)\to0$,
we have~$\dist(\Psi_*\kappa_{n_i},Q_k)\to0$, hence~$\dist(\Psi_*\mu,Q_k)=0$,
and~$\Psi_*\mu\in Q_k$ because~$Q_k$ is closed.
(d) Apply the argument of~(c) with~$\alpha_j:=(t,j)$; by~\eqref{eq:generators}, the
limit is~$t$.
\end{proof}

\begin{lemma}[Density] \label{lem:density}
The set~$\{\pi_\alpha\colon\alpha\in\CW\}$ is $L^1(\mu)$-dense in~$\Linf{\mu}$.
\end{lemma}

\begin{proof}
Let~$u\in\Linf{\mu}$ and~$\ve>0$.
By density of continuous functions in~$L^1(\mu)$, there exists a
continuous~$f\in C(\Omega)$ with~$\norm{u-f}_{L^1(\mu)}<\ve/2$.
By the Stone--Weierstrass theorem and approximation of real coefficients by
rationals, polynomials with rational coefficients in finitely many coordinates
are uniformly dense in~$C(\Omega)$.
Thus there are words~$\alpha_1,\ldots,\alpha_r \in \CW$
and a polynomial~$p\in\Q[z_1,\ldots,z_r]$ such that
\[
	\norm{f-p\circ (\pi_{\alpha_1},\ldots,\pi_{\alpha_r})}_\infty<\ve/2
	\qquad\text{on~$\Omega$}.
\]
In particular~$\norm{u-p\circ (\pi_{\alpha_1},\ldots,\pi_{\alpha_r})}_{L^1(\mu)}<\ve$.
Since~$\clip$ is $1$-Lipschitz and fixes the values of~$u$,
also~$\norm{u-\tilde p\circ (\pi_{\alpha_1},\ldots,\pi_{\alpha_r})}_{L^1(\mu)}<\ve$. By
\cref{lem:limit_law}(a), the left-hand side is~$\norm{u-\pi_{p(\alpha_1,\ldots,\alpha_r)}}_{L^1(\mu)}$.
\end{proof}

\begin{proof}[Proof of \cref{thm:existence}]
\emph{Definition of~$T$.}
We will first define~$T$ on $\{\pi_\alpha \colon \alpha\in\CW\}$,
then extend its definition to~$\Linf{\mu}$.
For~$\alpha\in\CW$, set~$T\pi_\alpha:=\pi_{T\alpha}$.
Note that if~$\pi_\alpha=\pi_\beta$ $\mu$-almost
everywhere, then~$\pi_{T\alpha}=\pi_{T\beta}$ $\mu$-almost everywhere,
by \cref{lem:limit_law}(b). By the same lemma, the map~$T$ so defined satisfies
\begin{equation} \label{eq:modulus_limit}
	\norm{Tu-Tw}_{L^1(\mu)}\leq\omega(\norm{u-w}_{L^1(\mu)})
\end{equation}
with~$u,w\in \{\pi_\alpha \colon \alpha\in\CW\}$, and it takes values
in~$\Linf{\mu}$.
By \cref{lem:density}, the map~$T$ extends uniquely to a continuous map on~$\Linf{\mu}$,
and such map satisfies~\eqref{eq:modulus_limit} for all~$u,w\in \Linf{\mu}$.

\emph{Convergence.} It remains to prove that~$T_n\to T$.
Fix~$k\geq1$. Recall that~$Q_k$ is the limit
of~$\overline{P_k(T_n)}$, thus it remains to prove~$\overline{P_k(T)}=Q_k$.

Let~$t\in Q'_k$ and~$u:=(\pi_{(t,1)},\ldots,\pi_{(t,k)})$.
Then~$Tu=(\pi_{T(t,1)},\ldots,\pi_{T(t,k)})$, so~$\CD_\mu(u,Tu)=t$ by
\cref{lem:limit_law}(d). Hence~$Q'_k\subseteq P_k(T)$, and taking closures
gives~$Q_k\subseteq\overline{P_k(T)}$.

Conversely, let~$u\in\Linf{\mu}^k$ and~$\ve>0$, and choose~$\delta>0$
with~$k(\delta+\omega(\delta))\leq\ve$. By \cref{lem:density}, there are
words~$\alpha_1,\ldots,\alpha_k$ with~$\norm{\pi_{\alpha_j}-u_j}_{L^1(\mu)}\leq\delta$ for
all~$1\leq j\leq k$. Let~$w:=(\pi_{\alpha_1},\ldots,\pi_{\alpha_k})$,
then~$Tw=(\pi_{T\alpha_1},\ldots,\pi_{T\alpha_k})$. By \cref{lem:elementary}(a)
and~\eqref{eq:modulus_limit},
\[
	\begin{split}
	d_W\bigl(\CD_\mu(u,Tu),\CD_\mu(w,Tw)\bigr)
	&\leq\sum_{j=1}^k\bigl(\norm{u_j-w_j}_{L^1(\mu)}+\norm{Tu_j-Tw_j}_{L^1(\mu)}\bigr)\\
	&\leq k\bigl(\delta+\omega(\delta)\bigr)\leq\ve,
	\end{split}
\]
and~$\CD_\mu(w,Tw)\in Q_k$ by \cref{lem:limit_law}(c). Hence~$\dist(\CD_\mu(u,Tu),Q_k)\leq\ve$
for every~$\ve>0$, and~$\CD_\mu(u,Tu)\in Q_k$ because~$Q_k$ is closed.
Thus~$P_k(T)\subseteq Q_k$, and~$\overline{P_k(T)}\subseteq Q_k$.
\end{proof}

\begin{remark}[Stars]\label{rem:star}
Let~$W_n$ average over neighbours on the star with centre~$0$ and leaves~$1,\ldots,n-1$,
and let~$\mu_n$ be uniform. The centre has vanishing weight but determines the output at
every leaf, so~$(W_n)$ is action convergent but not $L^1$-equicontinuous, and by
\cref{thm:existence} it has no limit in~$\CT$. A limit can still be defined using the
axiom of choice, via the Hahn--Banach theorem, but the resulting operator is not
continuous.
\end{remark}


\section{Continuous time}\label{sec:continuous_time}

The trajectory results extend to continuous time when~$F$ and~$G$ are Lipschitz.
Let~$W_n$ be an $L^1$-equicontinuous sequence with~$W_n\to W$,
as in \cref{thm:pairwise}, and let~$T_n=\Gamma(W_n)$ and~$T=\Gamma(W)$.
We consider the relaxation equations
\begin{equation}\label{eq:continuous_relaxation}
 \dot u_n=T_nu_n-u_n,\qquad \dot u=Tu-u,
\end{equation}
for curves~$t\mapsto u_n(t)$ in~$\Linf{\mu_n}$ and~$t\mapsto u(t)$
in~$\Linf{\mu}$.
The term~$-u$ makes the unit ball invariant, since~$T$ maps it into itself.
Systems without a bounded invariant set would require
comparing maps between balls of different radii; we do not pursue this.
For a matrix~$W_n$, the first equation is the system of coupled ODEs
\[
    \dot u(i,t)=F\p{u(i,t),\sum_jW_n(i,j)\,G(u(i,t),u(j,t))}-u(i,t),
\]
the continuous-time analogue of~\eqref{eq:intro_nonlin}.
With~$\tilde F(a,c):=F(a,c)-a$, this is the usual form
$\dot u(i,t)=\tilde F\p{u(i,t),\sum_jW_n(i,j)\,G(u(i,t),u(j,t))}$.
Conversely, if~$\tilde F$ is $K$-Lipschitz and points inwards,
that is~$\tilde F(1,c)\le0\le\tilde F(-1,c)$, then~$F(a,c):=a+\tilde F(a,c)/K$ takes
values in~$[-1,1]$, so after rescaling time by~$K$ every such system is of the
form above.

\subsubsection*{Path profiles}
Sampling~$x$ according to~$\mu$ and following~$t\mapsto u(t)(x)$ turns a solution
into a random path in~$C_\tau:=C([0,\tau],[-1,1])$, with the
metric~$d_\tau(a,b):=\sup_{t\le\tau}\abs{a(t)-b(t)}$. This space is Polish and bounded,
but not compact. On~$\CM(C_\tau)$ we use the Wasserstein distance~$d_W$,
with test functions that are $1$-Lipschitz for~$d_\tau$, and on its
subsets the resulting Hausdorff pseudometric~$d_H$.
For simplicity, we follow a single initial condition, the analogue of $1$-profiles.

A curve~$c\colon[0,\tau]\to\Linf{\mu}$, continuous for~$\norm{\cdot}_\infty$,
has a version with continuous paths: fix versions of~$c(q)$ for~$q$ in a countable
dense set~$D\subseteq[0,\tau]$; outside a single null set,
$\abs{c(q)(x)-c(q')(x)}\le\norm{c(q)-c(q')}_\infty$ for all~$q,q'\in D$, so~$q\mapsto c(q)(x)$
extends continuously to~$[0,\tau]$. Hence the law~$\CD_\mu(c)\in\CM(C_\tau)$ is well
defined, and~$d_W\bigl(\CD_\mu(c),\CD_\mu(c')\bigr)\le\sup_{t\le\tau}\norm{c(t)-c'(t)}_\infty$.

For a flow~$S=(S_t)_{t\ge0}$ on~$\Linf{\mu}$ with~$t\mapsto S_tu$ continuous
for~$\norm{\cdot}_\infty$, the \emph{path profile} of~$S$ on~$[0,\tau]$ is
\[
 P^{[0,\tau]}(S):=\bigl\{\CD_\mu\bigl((S_tu)_{0\le t\le\tau}\bigr)\colon
 u\in\Linf{\mu}\bigr\}\subseteq\CM(C_\tau).
\]
Evaluation at times~$0$ and~$\tau$ is $2$-Lipschitz from~$C_\tau$ to~$[-1,1]^2$ and
maps~$P^{[0,\tau]}(S)$ onto~$P_1(S_\tau)$.

\begin{proposition}[Continuous-time trajectories]\label{prop:continuous_time}
In the setting above, suppose in addition that~$F,G$ are Lipschitz.
The equations~\eqref{eq:continuous_relaxation} define unique global flows~$S^{(n)}$
and~$S$ on the unit balls, continuous in time for~$\norm{\cdot}_\infty$, and for
every~$\tau>0$,
\[
 d_H\bigl(P^{[0,\tau]}(S^{(n)}),P^{[0,\tau]}(S)\bigr)\longrightarrow0.
\]
Moreover, $S_\tau\in\CT$ and~$S^{(n)}_\tau\to S_\tau$ for every~$\tau>0$.
\end{proposition}

\begin{proof}
Let~$A$ be any of the maps~$T_n$ or~$T$, on~$\Linf{\mu}$, with~$\mu$ the corresponding measure.
By \cref{lem:linear_limit,lem:signed_kernel}, $A$ is given by the formula
of \cref{thm:pairwise}, with measures~$\nu_x$ such that~$\abs{\nu_x}(\Omega)\le1$. Given~$v,w\in\Linf{\mu}$, the set~$B$
where~$\abs{v-w}>\norm{v-w}_\infty$ is $\mu$-null, so~$\abs{\nu_x}(B)=0$ for
almost every~$x$, by \cref{lem:signed_kernel}. Hence~$A$ is $L$-Lipschitz
for~$\norm{\cdot}_\infty$, with~$L:=L_F(1+2L_G)$, where~$L_F,L_G$ are Lipschitz constants
of~$F,G$ for the sum metric.
The field~$A-\id$ is Lipschitz on the closed convex
set~$\Linf{\mu}\subseteq L^\infty$, and this set is invariant:
in variation-of-constants
form $v(t)=e^{-t}v(0)+\int_0^te^{s-t}Av(s)\d s$ is a convex combination of points
of~$\Linf{\mu}$. Picard iteration in~$C([0,\tau],\Linf{\mu})$ therefore gives, for all
times, a unique solution, and this solution stays in the ball. Since~$\norm{Av-v}_\infty\le2$,
each curve~$t\mapsto S_tu$ is $2$-Lipschitz in~$L^\infty$, so its path law is defined.

For~$h=\tau/m\le1$, the Euler
map~$E_h:=(1-h)\id+hA$ preserves~$\Linf{\mu}$. Let~$v_h$ interpolate linearly
between~$E_h^ju$ at times~$jh$. The classical Euler estimate (a local error~$2h$
propagated by Grönwall's inequality) gives
$\sup_{t\le\tau}\norm{S_tu-v_h(t)}_\infty\le Ch$ with~$C:=2(e^{(L+1)\tau}-1)$. The
constant~$C$ depends only on~$F$, $G$ and~$\tau$, not on~$n$, $u$ or~$h$.
Hence the path laws of~$(S_tu)_{t\le\tau}$ and~$v_h$ are within~$Ch$.

The map~$E_h$ is~$\Gamma$ with outer
function~$(a,c)\mapsto(1-h)a+hF(a,c)$. By \cref{thm:pairwise,lem:moduli},
the Euler maps of~$T_n$ form an $L^1$-equicontinuous sequence converging to that of~$T$,
which lies in~$\CT$. By \cref{thm:trajectories}, the sets of
laws~$\CD(u,E_hu,\ldots,E_h^mu)$, $u\in\Linf{\mu}$, converge in~$d_H$. Linear
interpolation from~$[-1,1]^{m+1}$ to~$C_\tau$ is $1$-Lipschitz, so the path profiles of
the interpolants converge too. Let~$\ve_n(h)\to0$ be their distance.
No uniformity in~$h$ is needed.

By the triangle inequality,
\[
 d_H\bigl(P^{[0,\tau]}(S^{(n)}),P^{[0,\tau]}(S)\bigr)\le2Ch+\ve_n(h).
\]
Let~$n\to\infty$, then~$h\to0$. The iterates~$E_h^m$ for~$A=T$ lie in~$\CT$ and
approximate~$S_\tau$ within~$Ch$ in~$L^1$, uniformly on the ball. A uniform limit of
$L^1$-uniformly continuous maps is $L^1$-uniformly continuous, so~$S_\tau\in\CT$.
Evaluating paths at times~$0$ and~$\tau$ gives~$P_1(S^{(n)}_\tau)\to P_1(S_\tau)$.
The same argument with~$k$ initial conditions, using \cref{thm:trajectories} for
general~$k$, gives the~$k$-profiles, so~$S^{(n)}_\tau\to S_\tau$.
\end{proof}

Thus, for Lipschitz~$F$ and~$G$, the examples of \cref{sec:applications} have
continuous-time limits.

\begin{remark}[Vlasov equation]\label{rem:continuous_vlasov}
In continuous time, the role of the self-consistent transfer operator
of \cref{subsec:STO} is played by the Vlasov equation
\[
 \partial_t\nu_t+\partial_z\bigl((\Phi_{\nu_t}-\id)\,\nu_t\bigr)=0,
\]
with~$\Phi_\nu$ as in~\eqref{eq:STO}. It describes the distribution~$\nu_t$ of a
population of particles, in which a particle at~$z$ moves with
velocity~$\Phi_{\nu_t}(z)-z$, so the velocity depends on the whole distribution.
Consider the system~$\dot u=T_nu-u$, with~$T_n$ as in~\eqref{eq:examples_dynamics}
and the degree measure. On the complete graph with loops,
the right-hand side is~$\Phi_\nu(u(i))-u(i)$, where~$\nu$ is the empirical law,
as in \cref{subsec:STO}. So
every vertex moves as a particle, and the empirical law solves the Vlasov equation
exactly, for every~$n$. On quasirandom graphs this holds only asymptotically. The
limit is~$J$ (\cref{subsec:quasirandom}), for which the same identity holds exactly,
and \cref{prop:continuous_time} transfers it to the finite systems, as \cref{cor:STO}
does in discrete time. The $\mu_n$-law of the solution therefore stays close to the
solution of the Vlasov equation with the same initial law, uniformly over initial
conditions and bounded time intervals.
\end{remark}



\section{Remarks, generalizations, and open problems}\label{sec:remarks}

\subsection{Lipschitz rates} \label{rem:rates}
The moduli of continuity in our proofs are obtained explicitly from the hypotheses.
If~$F$ and~$G$ are Lipschitz rather than merely continuous, that is, have linear
moduli~$\omega(s)=Ls$, the proofs give quantitative bounds. For instance, for coupled
maps of the form~$(T_nu)(i)=\sum_{j=1}^nW_n(i,j)\,h\bigl(u(j)\bigr)$, where~$h$ is
$L$-Lipschitz and~$\norm{W_n}_{1\to1}\leq L$ for all~$n$, the limit is~$Tu=W(h\circ u)$,
and from the proofs one extracts
\[
	d_H\bigl(P_1(T_n),P_1(T)\bigr)\leq(1+L+L^2)\,d_H\bigl(P_2(W_n),P_2(W)\bigr).
\]
This is the very special case~$F(a,b)=b$ and~$G(a,b)=h(b)$. For general~$G$, the
approximation by tensors in the proof of \cref{thm:pairwise} involves profiles of
increasing order, and we leave the corresponding quantitative bounds to future work.

\subsection{General state spaces}\label{rem:compact_states}
All systems in this paper take values in~$[-1,1]$, but the definitions and most
results extend to states in a compact metric space~$(X,\rho)$, such as~$[0,1]^d$ or
the circle~$\T=\R/\Z$ with its quotient metric. States are measurable
maps~$u\colon\Omega\to X$ modulo null sets, the~$L^1$ distance is replaced
by~$d_1(u,v):=\int\rho(u,v)\d\mu$, and profiles are subsets of~$\CM(X^{2k})$, compared
in the Wasserstein metric of~$X^{2k}$ with the sum metric. With these changes,
\cref{lem:elementary}, the repair lemma (\cref{lem:repair}) and all results of
\cref{sec:lifting} hold with the same proofs, with the constraint maps~$g_j$
allowed to go between different compact metric spaces.
This covers Kuramoto-type maps such
as~$\theta_i\mapsto\theta_i+\ve\sum_jW_n(i,j)\sin(2\pi(\theta_j-\theta_i))$. Indeed,
writing~$\sin(2\pi(\theta_j-\theta_i))=\sin(2\pi\theta_j)\cos(2\pi\theta_i)
-\cos(2\pi\theta_j)\sin(2\pi\theta_i)$, these maps are of the form in \cref{thm:several}
with~$N=2$, $G_1(\theta)=\sin(2\pi\theta)$, $G_2(\theta)=\cos(2\pi\theta)$,
and~$F(\theta,c_1,c_2)=\theta+\ve\bigl(c_1\cos(2\pi\theta)
-c_2\sin(2\pi\theta)\bigr)$ modulo~$1$.

\subsection{Relation with the linear theory} \label{rem:linear_relation}
The existence theorem of~\cite{backhausz2022action} concerns linear maps defined on the
unit ball~$\Linf{\mu}$ of~$L^\infty(\mu)$ with values in~$L^1(\mu)$, under a uniform
bound on~$\norm{\cdot}_{p\to q}$ for some~$p<\infty$. Here we require that maps send
the $L^\infty$-unit ball into itself, for the simple reason that such maps can be
iterated, and our focus is dynamics. In exchange, we only require
$L^1$-equicontinuity, which is implied by their bound: for linear maps of the
$L^\infty$-unit ball into itself, \cref{thm:existence} therefore covers the existence
theorem of~\cite{backhausz2022action}.
Combining the two proofs should give a common generalization, which we only sketch:
for~$p\in[1,\infty)$ and~$q\in[1,\infty]$, every action convergent sequence of
maps~$T_n\colon\Linf{\mu_n}\to L^q(\mu_n)$ with
\[
	\norm{T_nu-T_nw}_q\leq\omega(\norm{u-w}_p)
	\qquad\text{and}\qquad
	\sup_n\norm{T_n0}_q<\infty
\]
has a limit with the same modulus~$\omega$. As in~\cite{backhausz2022action}, compactness
is replaced by tightness and the Wasserstein metric by the L\'evy--Prokhorov metric.
We have not investigated the extension of \cref{thm:pairwise,thm:trajectories}
to this setting.

\subsection{Koopman operators}\label{rem:koopman}
Let~$f$ be a measure-preserving transformation of a standard probability
space~$(\Omega,\mu)$. Its Koopman operator~$T_fu:=u\circ f$ belongs to~$\CT$: it
maps~$\Linf{\mu}$ into itself and is an isometry of~$L^1(\mu)$. One may then
be tempted to study sequences of nonlinear maps through the action convergence of their
Koopman operators. This approach is not satisfactory, for two reasons. A modest one
is that the maps~$T\colon[-1,1]^n\to[-1,1]^n$ considered above do not come with a
natural invariant measure on~$[-1,1]^n$. A more fundamental one is that profiles
of Koopman operators are not very useful. One can check
that~$\overline{P_k(T_f)}=\overline{P_k(T_h)}$ for all~$k$ if and only
if~$f$ and~$h$ are weakly equivalent in the sense of
Kechris~\cite{kechris2010global}, that is, at distance zero in the partition metric
of Ab\'ert and Elek~\cite{abert2026space}. By the Rokhlin lemma, all
aperiodic measure preserving transformations are weakly equivalent. For instance,
profile closures cannot even distinguish an irrational rotation~$x\mapsto x+\alpha$
from the doubling map~$x\mapsto2x$ on~$\R/\Z$.
Without closures the equivalence is finer: $P_k(T_f)=P_k(T_h)$ for
all~$k$, or already for~$k=1$, if and only if~$f$ and~$h$ are factors of each
other, a property weaker than measurable isomorphism, but stronger
than weak equivalence.

\subsection{A soft regularity lemma} \label{rem:regularity}
By compactness of the nonempty closed subsets of $\CM([-1,1]^{2k})$ under~$d_H$, the
class~$\CT$ is totally bounded for the pseudometric
\[
	d(T,T'):=\sum_{k\geq1}2^{-k}d_H\bigl(P_k(T),P_k(T')\bigr).
\]
For every~$\ve>0$ there are finitely many~$T_1,\ldots,T_N\in\CT$ such that every~$T\in\CT$,
on any space, satisfies~$d(T,T_j)\leq\ve$ for some~$j$. Thus, up to~$\ve$, every coupled
system~\eqref{eq:intro_nonlin}, on a graph of any size and for any~$F,G$, behaves like one
of finitely many systems.
Making this quantitative, with an explicit bound on~$N$ in terms of~$\ve$ and an explicit
description of the maps~$T_j$, is left as an open problem.


\bibliographystyle{plain}
\bibliography{refs}

\end{document}